\documentclass[12pt]{amsart}
\usepackage{amsbsy}
\usepackage{graphicx,epsfig,subfigure,psfrag}
\begin{document}

%%%%%%%%%%%%%%%%%%%%%%%%%%%%%%%%%%%%%%%%%%%%%%%%%%%%%%%%%%%%%%%%%%%%
% Theorem, definition, lemma, proposition, corollary and proof
%%%%%%%%%%%%%%%%%%%%%%%%%%%%%%%%%%%%%%%%%%%%%%%%%%%%%%%%%%%%%%%%%%%%
%%%%%%%%%%%%%%%%%%%%%%%%%%%%%%%%%%%%%%%%%%%%%%%%
\newtheorem{theorem}{Theorem}%[section]
\newtheorem{proposition}{Proposition}%[section]
\newtheorem{lemma}{Lemma}%[section]
\newtheorem{corollary}{Corollary}%[section]
\newtheorem{definition}{Definition}%[section]
\newtheorem{remark}{Remark}%[section]
%%%%%%%%%%%%%%%%%%%%%%%%%%%%%%%%%%%
%%%%%%%%%%%%%%%%%%%%%%%%%%%%%%%%%%%%%%%%%%%%%%                  NEW
\newcommand{\beq}{\begin{equation}}
\newcommand{\eeq}{\end{equation}}
%%%%%%%%%%%%%%%%%%%%%%%%%%%%%%%%%%%%%%%%%%%%
\numberwithin{equation}{section}
\numberwithin{theorem}{section}
\numberwithin{proposition}{section}
\numberwithin{lemma}{section}
\numberwithin{corollary}{section}
\numberwithin{definition}{section}
\numberwithin{remark}{section}
%%%%%%%%%%%%%%%%%%%%%%%%%%%%%%%%%%%%%%%%%%%%
\newcommand{\re}{{\mathbb R}}
\newcommand{\n}{\nabla}
\newcommand{\ren}{{\mathbb R}^N}
\newcommand{\iy}{\infty}
\newcommand{\pa}{\partial}
\newcommand{\fp}{\noindent}
\newcommand{\ms}{\medskip\vskip-.1cm}
\newcommand{\mpb}{\medskip}
%%%%%%%%%%%%%%%%%%%%%%%%%%%%%%%%%%%%%%%%%%%%%%%%%
\newcommand{\BB}{{\bf B}}
\newcommand{\Am}{{\bf A}_{2m}}
%%%%%%%%%%%%%%%%%%%%%%%%%%%%%%%%%%%%%%%%%%%%%%%%%%%%
\renewcommand{\a}{\alpha}
\renewcommand{\b}{\beta}
\newcommand{\g}{\gamma}
\newcommand{\G}{\Gamma}
\renewcommand{\d}{\delta}
\newcommand{\D}{\Delta}
\newcommand{\e}{\varepsilon}
\renewcommand{\l}{\lambda}
\renewcommand{\o}{\omega}
\renewcommand{\O}{\Omega}
\newcommand{\s}{\sigma}
\renewcommand{\t}{\tau}
\renewcommand{\th}{\theta}
\newcommand{\z}{\zeta}
\newcommand{\wx}{\widetilde x}
\newcommand{\wt}{\widetilde t}
\newcommand{\noi}{\noindent}
%%%%%%%%%%%%%%%%%%%%%%%%%%%%%%%%%%

\title %%%[Blow-up localization]
{\bf Vast multiplicity of very singular  self-similar solutions of
a semilinear higher-order diffusion equation with
%%non-autonomous
time-dependent absorption}
 % with singular
%initial data
%%%%\footnote{\rm To appear in Commun. Contemp. Math.}}

\author{V.A. Galaktionov}
 %%% and A.E. Shishkov}

\address{Department of Mathematical Sciences, University of Bath,
 Bath BA2 7AY, UK}
 %%and Keldysh Institute of Applied Mathematics,
 %%%Miusskaya Sq. 4, 125047 Moscow, RUSSIA}

\email{vag@maths.bath.ac.uk}

%%\address{Institute of Applied Mathematics and Mechanics
%%of NAS of Ukraine, R. Luxemburg str. 74, 83114 Donetsk,
%%UKRAINE}

%%%\email{shishkov@iamm.ac.donetsk.ua}

%%%%\thanks{Research partially supported by the INTAS Project 03-51-5007}
%\\ \quad
%{To appear in Commun. Contemp. Math.}}

 \keywords{The Cauchy problem,
diffusion equations with absorption, initial Dirac mass, very
singular solutions, existence, nonexistence, bifurcations, branching. }
 \subjclass{35K55, 35K40, 35K65. }
\date{\today}

%%%%%%%%%%%%%%%%%%%%%%%%%%%%%%%%%
% \newcommand{\begin{equation}}{\begin{equation}}
%%%%%%%%%%%%%%%%%%%%%%%%%%%%%%%%%%%%%%%%%%%%%%%
% \newcommand{\end{equation}}{\end{equation}}
%%%%%%%%%%%%%%%%%%%%%%%%%%%%%%%%%%%%%%%%%%%%%%%%%%%%%%%%%%%

%%%%%%%%%%%%%%%%%%%%%%%%%%%%%%%%%%%%%%%%%%%%%%%%%%%%%%%%

% \vskip 1.5cm

\begin{abstract}
As a basic model, the Cauchy problem in $\ren \times \re_+$ for
the $2m$th-order semilinear parabolic equation of the
diffusion-absorption type
%%%% with non-autonomous sink-term,
 %of the form
 $$
u_t =-(- \D)^m u - t^\a |u|^{p-1}u, \quad \mbox{with}
 \,\, p > 1, \,\, \a > 0, \,\,\, m \ge 2,
 $$
 %%% is studied, and
%%  the Cauchy problem in $\ren \times \re_+$
 with
 {\em singular} initial data $u_0 \not =0$ such that
 $
u_0(x)=0$ for any $x \not = 0$ is studied.
  The additional multiplier
  %%% $h(t)=t^\a$ in the absorption term
 %%% satisfies
$h(t)=t^\a \to 0$ as $t \to 0$ in the absorption term plays a role
of time-dependent non-homogeneous potential that affects the
strength of the absorption term in the PDE. Existence and
nonexistence of the corresponding very singular solutions (VSSs)
is studied.
%%% of the Cauchy problem.
For $m=1$ and $h(t) \equiv 1$, first nonexistence result for $p
\ge p_0=1+\frac 2N$ was proved in the celebrated paper by Brezis
and Friedman in 1983. Existence of VSSs in the complement interval
$1<p<p_0$ was established  in the middle of the 1980s.

 %% In particular,
 %% Our main goal is to  justify that,
   %%, for the power case $h(t)=t^\a$ with $\a>0$
 The main goal is to justify that,  in the subcritical range $1<p<p_0=1+\frac{2m(1+\a)}N$, there
 exists
 a finite number
  of different VSSs of
  the self-similar form
 $$
 \mbox{$
 u_*(x,t) = t^{-\b} V(y), \quad y = x/t^{\frac 1{2m}}, \quad \b= \frac{1+\a}{p-1},
 $}
 $$
 where each $V$ is an exponentially decaying as $y \to \infty$ solution of the
elliptic equation
  $$
  -(-\D)^m V + \mbox{$\frac 1{2m}$} \,y \cdot  \nabla V  + \b V -
  |V|^{p-1}V=0 \quad \mbox{in} \quad \ren.
  $$
  Complicated families of VSSs in 1D and
  also  non-radial VSS patterns in $\ren$ are detected. Some of these
  VSS profiles $V_l$ are shown to bifurcate from $0$ at the
 bifurcation exponents
 $$
 \mbox{$
 p_l = 1+\frac{2m(1+\a)}{l+N}, \quad \mbox{where} \,\,\,l=0,1,2,...
 \,.
  $}
$$
%%%% are studied.
%% VSS profiles are also convenient to detect at
%% $\a$-bifurcation points.
 %%%The limit $\a \to +\infty$ is also considered
 %%and describe the final ($\a=+\infty$) profile.
%%% which allows to predict some properties of solutions with
 %%we naturally arrive at
 %%%non-power potentials such as
%% $
%% h_*(t)= {\mathrm e}^{-1/t},
 %%%  ???? ??? ??? ???
 %% $
 %%  exhibiting  new evolution
%%%
 %%%% properties.
   %%%????
 %%% where
 %%%the first one $V_0$ is shown to be stable. Discrete and countable
 %%subsets of other self-similar and approximately self-similar
  %%%patterns are introduced.

\end{abstract}

%%%%%%%%%%%%%%%%%%%%%%%%%%%
\maketitle

%\vskip 0.7cm

%%%%%%%%%%%%%%%%%%%%%%%%%%%%%%%%%%%%%%%%%%%%%%%%%%%%%%%%%%%%%%%%%%%%%%%%
% INTRODUCTION
%%%%%%%%%%%%%%%%%%%%%%%%%%%%%%%%%%%%%%%%%%%%%%%%%%%%%%%%%%%%%%%%%%%%%%%%
\setcounter{equation}{0}
\section{Introduction: VSSs for
%%%Brezis--Friedman results  and extensions to
higher-order diffusion-absorption PDEs}
 \label{Sect1}
  \setcounter{equation}{0}

\subsection{Diffusion-absorption model}

Our basic model is
%%the Cauchy problem for
 the $2m$th-order
semilinear heat equation in $\ren \times \re_+$ with
non-autonomous (non-homogeneous) absorption term
 \beq
 \label{1.1R}
 u_t = -(-\D)^m u - h(t) |u|^{p-1} u \quad (p>1,  \,\,\,m \ge 2),
 \eeq
where the given function $h(t)$ satisfies
 $$
h(t) > 0, \quad h(t) \to 0 \quad \mbox{as \,$t \to 0$}.
 $$

 For
$m=1$ and $h=1$, the critical Fujita-like exponent $p_0=1 + \frac
2N$  and existence-nonexistence of {\em very singular solutions}
(VSSs)
 was derived since the 1980s; see first principal results of the 1980s in
 \cite{BrFr, GKS, BPT, EscKav, KPApp, KPe, KVer}, and a long  list of
  references in \cite[Ch.~4]{AMGV}.
In the higher-order case $m \ge 2$ and $h =1$, such an
investigation has been performed in \cite{Gal2m, GSVSS1, GW2}. The
influence to VSSs of non-homogeneous potentials $h(t)$ was first
studied in  Marcus--V\'eron \cite{MV02} in the case
%%in the second-order case
$m=1$; see more recent extensions in \cite{SV1}.

%%%%%%%%%%%%%%%%%%%%%%%%%%%%%%%%%%%%%%%%
\subsection{Similarity VSSs: towards a non-variational problem}

The layout of the paper is as follows. We  investigate the
phenomena associated with the non-homogeneous term and consider
the power case
 \beq
\label{al.1}
 \mbox{$
 h(t) = t^\a, \quad \mbox{with the fixed exponent}
\,\,\, \a>-1 \quad\big(\mbox{so that} \quad \int_0 h(t)\,{\mathrm
d}t < \infty\big).
 $}
 \eeq
%% The last condition guarantees the necessary integrability,
%%% $\int_0 h(t)\,{\mathrm d}t < \infty$.
%%%% In Section \ref{Sect2} we study similarity VSSs in the power case
 In Section \ref{Sect2},
  %%% we study similarity VSSs in the power case
for the higher-order equation (\ref{1.1R}) with $m \ge 2$, we
study the existence
 and multiplicity of similarity solutions and show
that, in the subcritical range $p \in (1,p_0)$, with the critical
(Fujita-like) exponent
 $$
 \mbox{$
 p_0= 1+ \frac{2m(1+\a)}N,
  $}
  $$
there exist VSSs
of the similarity form  %%%%%%(cf. (\ref{1.5R}))
 \beq
 \label{1.6R}
 \mbox{$
 u_*(x,t) = t^{-\b} V(y), \quad y = x/t^{\frac 1{2m}}, \quad \mbox{with} \,\,\, \b= \frac{1+\a}{p-1}>0,
  $}
 \eeq
 where $V$ is a non-trivial radial solution of the elliptic equation
%%\begin{align}
 \beq
\label{1.7R}
 \left\{
%%{\bf P}_{2m}(V) \equiv
 \begin{matrix}
\BB_1 V -  |V|^{p-1}V \equiv  -(-\D)^m V + \mbox{$\frac 1{2m}$}
\,y \cdot  \nabla V  + \b V - |V|^{p-1}V=0 \quad \mbox{in} \quad
\ren, \smallskip\smallskip
\\
 V(y) \,\,\, \mbox{decays exponentially fast as} \,\,\, |y| \to
\infty. \qquad\qquad\qquad\qquad\qquad\qquad\quad
 \end{matrix}
 \right.
 \eeq
 %%%\qquad\qquad
 %%%\label{ExpC}
%%\end{align}
 Passing to the limit $t \to 0^+$ in (\ref{1.6R}) yields  the initial function
 satisfying
 \beq
   \label{g15}
u(x,0)=0 \quad \mbox{for all} \,\,\, x \not = 0,
 \eeq
 which nevertheless is not a finite measure, and this reinforces the
 notion of {\em very singular solutions} for (\ref{1.6R}).
%%% in the sense of distributions,
 %%   in the sense of
%%%$(\ref{we1})$,
Condition in (\ref{1.7R}) on exponential decay at infinity is
naturally enforced by introducing weighted $L^2$ and Sobolev
spaces; details are given in  in  Section \ref{Sect2}.

%%%Note that, for any $m \ge 2$ the problem (\ref{1.7R}) is not
%%variational, the principle linear  operator ${\bf B}_1$ is not
%%self adjoint, and represents no order-preserving properties of the
%%corresponding parabolic flow in view of the lack of the Maximum
%%Principle (applies for $m=1$ only). These negative properties
%%dictates the specific structure of our analysis.

%%%%%%%%%%%%%%%%%%%%%%%%%%%%%%%%%%%%%%%%%
\subsection{On main results}
 %%%% approaches, and style of analysis}

 The main goal  is to detect
the whole range of VSS profiles by
 using spectral properties of the linear operator ${\bf B}_1$ in (\ref{1.7R}) and
 bifurcation theory. This determines a countable number of
 $p$-bifurcation branches of VSS profiles in the subcritical range
 $p<p_0$ that, for any $\a>-1$, are originated at the bifurcation points
  \beq
  \label{bif1}
  \mbox{$
  p_l = 1+ \frac{2m(1+\a)}{N+l}, \quad l=0,1,2,... \, .
   $}
    \eeq
This also yields   $\a$-bifurcation points,
%%% that, according to
 %%the linear operator
 %%%%(\ref{BB*}), are
 \beq
  \label{al11}
  \mbox{$
  \a_l = \frac{(p-1)(N+l)}{2m}-1, \quad l=0,1,2,... \, .
   $}
    \eeq

 %% analytical
 %%  and numerical methods that in the {radial setting},
 %%the ODE problem (\ref{1.7R}), (\ref{ExpC}) admits at least
 %%\beq
 %%\label{M**}
 %%\mbox{$
 %%M(m,p,N,\a)= \sharp_{\rm even}(  \mbox{$\frac{N(p_0-p)}{p-1}$}
%% )
%%  $}
%% \eeq
%%different non-trivial radial solution $f= f(|y|)$, where $
%%\sharp_{\rm even}(z)$ for $z \ge 0$ denotes the number of
%%non-negative even
%% numbers $0,2,4,...$ not exceeding the integer part  $\lfloor z \rfloor$.

%%  We also study the asymptotic stability of the VSS and describe a
%  countable subset of other self-similar or approximately
% self-similar patterns in the Cauchy problem (\ref{1.1R}).
%%NEW
%We show that in the supercritical range $p>p_0$, no generically
%stable (see a precise definition in Section \ref{Sect4})
%non-trivial VSSs solutions exist. On the other hand, there exists
%an uncountable family of different similarity solutions which are
%not in $L^1$ and have special stability properties.

%%For
%%instance, for $p=2$ and $\a \sim 20-40$, there exist hundreds
%% This
%%phenomenon of splitting is not available for $\a=0$ studied
%%earlier in \cite{GW2}.

In Section \ref{SectAl}, using numerical experiments, we show
that, surprisingly, the whole picture of branches of VSS
similarity profiles is essentially
 more complicated than that for $\a=0$ studied in \cite{GW2}. It
turns out that even in 1D for sufficiently large $\a>0$, the whole
set of VSS profiles consists of various families of solutions of
different geometric shapes. In addition, this
 includes new
branches that are expected to appear at saddle-node (turning)
bifurcation points in $\a$.
%% obtained via $\a$-bifurcations
In Section \ref{SectInf}, for future use for non-power potentials
such as
 \beq
  \label{al.2}
   h_*(t) = {\mathrm e}^{- 1/t} \, ,
 \eeq
 we are particularly interested in the study of the limit behaviour of such VSSs as $\a \to +\infty$.
  Some nonexistence results of VSSs  in the range $p \ge p_0$ are
proved in Section \ref{SectNon}.

%%%Section \ref{Sect3} is devoted to non-power potentials with typical behaviour
%%% such as

%%We are going to present a sharp estimate on $h(t)$ that guarantee
%%existence and non-existence of VSSs. The function $h_*(t)$ in
%%(\ref{al.2}) precisely corresponds to the critical case. ???? We
%%describe the structure of a ``boundary layer", which occurs at
%%$t=0$. ??? ???

\smallskip

It is necessary to mention that for, any $m \ge 2$, the problem
(\ref{1.7R}) even in 1D is not variational (it is for $m=1$ only),
so the operators there are not potential in any topology.
Moreover, even the linear principal part ${\bf B}_1$ therein in
not self-adjoint (symmetric) in any weighted space $L^2_\rho$. In
addition, the corresponding parabolic flow  is not
order-preserving  in view of the lack of the Maximum Principle
(again applies for $m=1$ only).

\smallskip

 Therefore, for describing the complicated discrete sets
$\{V_{\s_k}\}$, where $\s_k$ is a multiindex (see Section
\ref{SectAl}), we are going to use a machinery of various
bifurcation-branching methods of nonlinear analysis.
 Nevertheless, many of our final conclusions on the global behaviour of bifurcation
diagrams of $V$'s remain formal, and we do not expect that a
reasonably simple justification can be achieved soon. Therefore,
we heavily rely on numerical construction of, at least,  am couple
of hundreds  of similarity profiles. These numerics are also not
that easy at all.
%%%%though are often the
Often, the reliable numerics will be the only tool to detect this
{\em vast multiplicity} of VSS profiles of (\ref{1.7R}) far away
from bifurcation points.

%%%%%%%%%%%%%%%%%%%%%%%%%%%%%%%%%%%%%%%%%%%%%%%%%%%%%%%%%%%%
%%%%%%%%%%%%%%%%%%%%%%%%%%%%%%%%%%%%%%%%%%%%%%%%%%%%%%%%%%%%%%%%%
\section{VSSs profiles: local bifurcation theory}
\label{Sect2}

Thus,
 consider the Cauchy problems for the PDE (\ref{1.1R}) with
the power potential (\ref{1.1R}).
%% \beq
%% \label{1.1Rh}
%% u_t = -(-\D)^m u - t^\a |u|^{p-1} u, \quad u_0 \in L^1(\ren) \cap L^\infty(\ren).
%% \eeq

%%%%%%%%%%%%%%%%%%%%%%%%%%%%%%%%%%%%%%%%%%%%%%%%%
\subsection{The fundamental solution of the poly-harmonic equation}

%% We will show that  in the supercritical range $p>p_0$,
%%for a class of sufficiently small  initial data, the solutions of
%%(\ref{1.1R}) behave, as $t \to \infty$, as the

%%%This is
%% Consider the  fundamental solution
%%(up to a constant multiplier $C
%%%\not = 0$ specified by initial data),

 This is
  \beq
  \label{1.3R}
   b(x,t) = t^{-\frac N{2m}}F(y), \quad y= x/t^{\frac 1{2m}},
 \eeq
  being the  {\em fundamental solution}
 of the linear {\em poly-harmonic equation}
  \beq
  \label{Lineq}
  u_t = - (-\D)^m u.
  \eeq
  %%\\ VAG: we will use this as equation for $u$, not $b$
 The rescaled kernel $F$ is then the unique radial solution of the elliptic equation
  \begin{equation}
\label{ODEf}
 \mbox{$
  {\bf B} F \equiv -(-\Delta )^m F + {\mathcal L}_0 F = 0
 \,\,\,\mbox{in} \,\, \ren,  \,\,\, \int F = 1; \quad
 {\mathcal L}_0=  \mbox{$\frac 1{2m}$} \,y \cdot
\nabla  + \mbox{$\frac N{2m}$} \,I,
 $}
\end{equation}
satisfying for some positive constants $D>1$ and $d>0$  depending
on $m$ and $N$
 \cite{EidSys}
\begin{equation}
\label{fbar}
 |F(y)| < D \bar F(y) \equiv D \o_1 {\mathrm e}^{-d|y|^{\alpha}}
\,\,{\rm in} \,\, \ren, \quad \a=  \mbox{$\frac {2m}{2m-1}$} \in
(1,2),
\end{equation}
 $\o_1>0$ being a normalization constant such that $\int \bar F=1$.

It turns out that
 the linear operator $\BB_1$ in equation (\ref{1.7R}) is connected
with the operator (\ref{ODEf}) for the rescaled kernel $F$ as
follows:
%% of the poly-harmonic operator
 %%%(to be introduced in Section \ref{Sect2})
  %% (see
 %%%(\ref{1.3R}))
 \beq
 \label{BB*}
  \BB_1 = \BB + c_1 I, \quad \mbox{where} \,\,\, c_1 =
  \mbox{$\frac {N(p_0-p)}{2m(p-1)}$}.
  \eeq

%%%%%%%%%%%%%%%%%%%%%%%%%%%%%%%%%%%%%%%%%%%%%%%%%%%%%%%%%%%%%%%%%%%%%%%%%%%%%
\subsection{\bf The point spectrum of  the non self-adjoint
operator ${\bf B}$}

 %%For any  $m>1$, ${\mathbf B}$  is not symmetric and does not admit
 %%a self-adjoint extension.
  We consider the linear operator $\BB$ given in (\ref{ODEf}) in the weighted space
$L^2_\rho({\bf R}^N)$ with the exponentially growing weight
function
 \beq
  \label{rho44}
  \rho(y) = {\mathrm e}^{a |y|^\a}>0 \quad {\rm in} \,\,\, \ren,
 \eeq
  where $a \in (0, 2 d)$ is a sufficiently small
constant. We ascribe to ${\bf B}$ the domain $H^{2m}_\rho({\bf
R}^N)$ being
%%\\ Some change
%%NEW
 a Hilbert space with the
norm
\[
\mbox{$
\|v\|^2 = \int \rho(y) \mbox{$\sum_{k=0}^{2m}$} \, |D^{k}
 v(y)|^2 \, {\mathrm d} y,
  $}
\]
induced by the corresponding inner product. We have $H^{2m}_{\rho}
\subset L^2_{\rho}  \subset L^2 $.  The spectral properties ${\bf
B}$
%%NEW
are as follows \cite{Eg4}:

\begin{lemma}
\label{lemspec} {\rm (i)} ${\mathbf B}: H^{2m}_\rho \to L^2_\rho$
is a bounded linear operator with the real point spectrum
\begin{equation}
\label{spec1}
 \mbox{$
  \sigma({\mathbf B}) = \big\{\lambda_l = -\mbox{$\frac
l{2m}$}, \,\, l = 0,1,2,...\big\}.
 $}
\end{equation}
The eigenvalues $\l_l$ have finite multiplicity with eigenfunctions
\begin{equation}
\label{eigen}
 \mbox{$
 \psi_\beta(y) = \frac{(-1)^{|\b|}}{\sqrt {\b !}} \,
D^\beta F(y), \quad \mbox{with any} \,\,\,|\b|=l\ge 0. $}
\end{equation}

{\rm (ii)} The set of eigenfunctions $\Phi = \{\psi_\b, \, |\b| =
0,1,2,...\}$ is complete and closed
 in $L^2_\rho$.
\end{lemma}

 In the classical second-order case
$m=1$, we have that
 $$
 F(y) = (4 \pi)^{-\frac N2} {\mathrm e}^{-\frac 14 \, {|y|^2}}
  $$
   is the
 rescaled positive Gaussian kernel and
the eigenfunctions are
 $$
 \psi_{\beta}(y)= {\mathrm e}^{-\frac 14 \,|y|^2} H_\b(y),
  $$
  where $H_\b$ are separable
Hermite polynomials in ${\bf R}^N$ \cite[p.~48]{BS}. The operator
${\bf B}$, with the domain $H_\rho^2$ and the weight $\rho =
{\mathrm e}^{|y|^2/4}$, is self-adjoint and the eigenfunctions
form an orthogonal basis in $L^2_\rho$.

Lemma \ref{lemspec} gives  the   centre and stable subspaces of
$\BB$,
$
E^c = {\rm Span}\{\psi_0= F\}$, $E^s = {\rm Span}\{\psi_\b, \,
|\b|>0\}$.

%%%%%%%%%%%%%%%%%%%%%%%%%%%%%%%%%%%%%%%%
\subsection{\bf The polynomial eigenfunctions of the
adjoint operator ${\bf B}^*$}

Consider the adjoint operator to ${\bf B}$,
\begin{equation}
\label{B2}
 \mbox{$
 {\bf B}^* = -(-\D)^m - \mbox{$\frac 1{2m}$} \, y \cdot
\nabla \, .
 $}
\end{equation}
For $m=1$, the following representation holds:
%\begin{equation}
%\label{Bsymm}
 $$
  \mbox{$
  {\bf B}^* \equiv \frac 1{\rho^*} \, \nabla \cdot (\rho^* \nabla
), \,\,\, \mbox{with} \,\,\, {\mathcal D}({\bf B}^*) =
H^2_{\rho^*} , \quad \mbox{where} \quad
 %%% with weight $
\rho^*(y) = {\mathrm e}^{- \frac 14 \,|y|^2},
  $}
 $$
%\end{equation}
 is self-adjoint in  $L^2_{\rho^*} $   and has a discrete
 spectrum. The eigenfunctions form an orthonormal basis
in $L^2_{\rho^*} $ and the classical Hilbert-Schmidt theory
applies~\cite{BS}.

For $m>1$, we  consider $\BB^*$ in $L^2_{\rho^*} $ with the
exponentially decaying weight function
 $$
  \mbox{$
   \rho^*(y) = \frac
1{\rho(y)} \equiv {\mathrm e}^{-a|y|^\a} > 0.
 $}
 $$

\begin{lemma}
\label{lemSpec2} {\rm (i)}  $\BB^*: H^{2m}_{\rho^*}  \to
L^2_{\rho^*} $ is a bounded linear operator with the same spectrum
as $\BB$, $(\ref{spec1})$. The eigenfunctions $\psi_\b^*(y)$ with
$|\b|=l$ are $l$th-order polynomials
\begin{equation}
\label{psidec}
 \mbox{$\psi_\b^*(y) =  \frac {1}{ \sqrt{\b !}}\big[ y^\b +
\mbox{$\sum_{j=1}^{\lfloor \frac {|\b|}{2m} \rfloor} \frac 1{j
!}$}(-\Delta)^{m j} y^\b \big].
 $}
\end{equation}

{\rm (ii)} The set $\{\psi_\beta^*\}$ is complete in $L^2_{\rho^*}
$.
\end{lemma}

From this definition of the adjoint eigenfunctions, the
orthonormality condition holds
 \beq
 \label{Ortog}
\langle \psi_\b, \psi_\g^* \rangle = \d_{\b \g},
 \eeq
where $\langle \cdot, \cdot \rangle$ denotes the standard $L^2$
inner product.  For $m=1$, both (\ref{spec1}) and (\ref{psidec})
are well-known properties
 of the separable  Hermite polynomials,
  %% generated by a self-adjoint
%%Sturm-Liouville problem
\cite{BS}.

%%%%%%%%%%%%%%%%%%%%%%%%%%%%%%%%%%%%%%%%%%%%%
\subsection{Bifurcations at $p = p_l$: local existence of the VSSs}

It follows from (\ref{BB*}) and the above lemmas that  the only possible bifurcation points
of VSSs are obtained from
%%%%%the equality
 \beq
 \label{pcrit}
  \mbox{$
 c_1=- \l_l \quad \Longrightarrow \quad p=p_l= 1+ \frac{2m(1+\a)}{N+l}, \,\,\ l =0,1,2,... \, .
  $}
  \eeq

We justify these bifurcation phenomena. Taking $p$ near the
critical values as defined in (\ref{pcrit}), we look for small
solutions of the problem (\ref{1.7R}). At $p = p_l$, the linear
operator $\BB_1$ has a nontrivial kernel, hence, the following
result:

\begin{proposition}
% Let $N < 2m$.
 Let for an integer $l \ge 0$, the
eigenvalue $\l_l = -\frac l{2m}$ of operator $(\ref{ODEf})$ be of
 odd multiplicity. Then the exponent in $(\ref{pcrit})$
   is a bifurcation point for
 $(\ref{1.7R})$.
%%%%, $(\ref{ExpC})$.
\end{proposition}

\noindent {\em Proof.} This result is standard in bifurcation
theory; see \cite[p.~381]{Deim} and similar results in a related
problem in \cite[\S~6]{GW2}.
 Necessary spectral properties of
the linearized operator are given in \cite{Eg4}.
 We present some
comments concerning the linearized operators involved.

%% For a moment, given an $n \gg 1$, we denote
%%by $(|V|^{p-1}V)_n$ a suitable uniformly Lipschitz continuous
%%truncation of the nonlinearity $|V|^{p-1}V$ such that $
%% (|V|^{p-1}V)_n \equiv |V|^{p-1}V$ for $|V| \le n$ so that
%% $$
%% (|V|^{p-1}V)_n \to |V|^{p-1}V \quad \mbox{as} \,\,\, n \to \infty \,\,\,
%% \mbox{uniformly on compact subsets}.
%% $$

 Consider
 in $H^{2m}_\rho $ the equation
  \beq
  \label{treq1}
 \hat{\BB} V = -(1+c_1)V + |V|^{p-1}V,\quad \hat{\BB} =
\BB_1 - (1+c_1)I \equiv \BB-I.
 \eeq
The spectrum of $\hat{\BB}$ is a translation of that of $\BB$,
 %%(\ref{Sp1}),
  $\s(\hat{\BB})=\{-1- \frac l{2m}\}$, and consists of
strictly negative eigenvalues. The inverse operator
$\hat{\BB}^{-1}$ is known to be compact (Proposition 2.4 in
\cite{Eg4}). Therefore, in
 the corresponding integral equation
 \beq
 \label{Inteq}
 V = {\bf A}(V) \equiv - (1+c_1) \hat {\BB}^{-1} V +
\hat{\BB}^{-1} |V|^{p-1}V,
 \eeq
%%the right-hand side contains a compact Hammerstein operator,
%%%\cite[Ch.~V]{Kras}  (see also some details in \cite{BGW1}).
the right-hand side contains a compact Hammerstein operator in
$L^q_\rho(\ren)$ space for some $q  \ge 1$ \cite[p.~38]{Kras} (see
details on the resolvent of ${\bf B}$  in \cite{BGW1}). In this
application,  there exist certain technical difficulties in
checking compactness of this Hammerstein operators in weighted
$L^q$-spaces over whole $\ren$.  As an alternative, we can use
 Ladyzhenskii's theorem \cite[p.~34]{Kras} establishing
 compactness in $C$.
 %%% (in this case, the auxiliary truncation of the
 %%%nonlinearity is not necessary).

To avoid this technicalities, we use \cite[Thm.~28.1]{Deim} where
no assumptions on compactness of the vector field are necessary.
The main hypothesis therein  is oriented to the linearized
operator in (\ref{Inteq}) that is assumed to be Fredholm of index
zero at bifurcation values (which is true by \cite{Eg4}).  Then
the result is first  obtained for truncated uniformly bounded
nonlinearity instead of $|V|^{p-1}V$ (see \cite[p.~1088]{GW2}),
for which all the hypotheses are valid.
 Then
bifurcations in the truncated problem (\ref{Inteq}) are always
guaranteed if the derivative ${\bf A}'(0)=- (1+c_1)
\hat{\BB}^{-1}$ has the eigenvalue $1$ of odd multiplicity.
%% see
%%%\cite[???]{KrasZ} and
(cf. also \cite[p.~196]{Kras} for compact integral operators).
 %% For exponentially decaying kernels
 %% of $\BB^{-1}$, the case of unbounded space $\ren$ can be settled
 %% by approximation by a converging in the norm sequence of compact
%%operators in expanding bounded domains. Some of the compactness
%%conditions in \cite[Ch.~1]{Kras} are valid for arbitrary unbounded
%%domains in the integral operators.

Thus,  bifurcations in the  problem (\ref{Inteq}) occur if the
derivative ${\bf A}'(0)=- (1+c_1) \hat{\BB}^{-1}$ has the
eigenvalue $1$ of odd multiplicity.
 %% see \cite{KrasZ} and
%%\cite{Kras}.
 Since
 $
 %% \mbox{$
  \s({\bf A}'(0))= \{(1+c_1)/(1+ \frac
l{2m})\},
 %%% $}
 $
  we arrive at critical values (\ref{pcrit}). By
construction, the solutions of (\ref{Inteq}) for $p \approx p_l$
are small in $L^2_\rho $ and, as can be seen from the properties
of the inverse operator, in $H^{2m}_\rho$.
 Since the weight $\rho(y)$ is a monotone growing function as $|y| \to
\infty$, this implies that  $V \in H^{2m}_\rho $ in the
subcritical Sobolev range
 \beq
 \label{sob}
  \mbox{$
 1<p<p_S= \frac{N+2m}{N-2m}
 $}
 \eeq
 (which is true in our VSS range $p \le p_0$)
 is a
uniformly bounded, continuous function by standard elliptic
regularity results and embedding theorems;  these can be found in
\cite{Maz} and \cite{Tay}. For $N<2m$, the result is
straightforward in view of embedding \cite[p.~5]{Tay}
 $$
  \mbox{$
 H^m(\ren) \subset C(\ren) \quad (N<2m).
 $}
 $$
 Let us mention
related boundedness results of  parabolic orbits of (\ref{1.7R})
($\a=0$) in the range (\ref{sob})  obtained  in \cite[\S~2]{GW2}
by Henry's version of Gronwall's inequality with power kernel and
scaling arguments.

%% (as we have mentioned, for compactness in $C$, this is
%%%not necessary).
 %%%Since the weight (\ref{rho44}) is a monotone growing function as $|y| \to
%%%\infty$, using the known asymptotic properties of solutions of the
%%ODE (\ref{1.7R}) (Section \ref{Sect2}), $V \in H^{2m}_\rho $ is a
%%uniformly bounded, continuous function. This well corresponds to
%%%local results of classic elliptic theory.
%%[It is worth mentioning
%%that, for even $m$, solutions of (\ref{1.7R}) may blow-up at
%%finite $y$ (a striking contrast to second-order ODEs) forming
%%singularities $\not \in L^2_\rho$ locally.]

 Therefore, for $p \approx p_l$, we have bounded, uniformly small solutions of (\ref{Inteq}) only.
By interior regularity results for elliptic equations,  these
solutions are smooth enough to be classical ones of the
differential equation (\ref{treq1}).
 %% Hence the same
%%bifurcations occur in the original non-truncated equation
%%(\ref{Inteq}) corresponding to $n=\infty$.
    $\qed$

\smallskip

 Thus, $l=0$ is always a bifurcation point since $p_0=0$ is simple.
 In general, for $l=1,2,...$, the odd multiplicity occurs depending on the
 dimension $N$.
In particular, for $l=1$, the multiplicity is $N$, and for $l=2$,
it is $\frac{N(N+1)}2$.
  In the case of the even multiplicity of
 $\l_l$, an extra analysis is necessary to guarantee that a
 bifurcation occurs, \cite{KrasZ}. It is important that, for
 key applications, namely, for $N=1$ and for the radial setting in $\ren$,
 the eigenvalues
(\ref{spec1}) are simple and (\ref{pcrit}) are  bifurcation
points.

Since the nonlinear perturbation term in the integral equation
(\ref{Inteq}) is an odd sufficiently smooth operator, we easily
obtain the following  result describing the local behaviour of
bifurcation branches, see  \cite{Kras} and \cite[Ch.~8]{KrasZ}.

\begin{proposition}
\label{PBif2}
 Let $\l_l$ be a simple eigenvalue of $\BB$ with  eigenfunction $\psi_l$.
  Denoting
 \beq \label{kap1}
\kappa_l = \langle |\psi_l|^{p-1}\psi_l, \psi_l^* \rangle,
 \eeq
  we have that
problem $(\ref{1.7R})$
%%%, $(\ref{ExpC})$
 has
 {\rm (i)} precisely two small solutions  for $p \approx p_l^-$ and
no solutions for $p \approx p_l^+$
 if $\kappa_l >0$, and {\rm (ii)} precisely two small solutions
 for $p \approx p_l^+$ and no solutions for $p \approx p_l^-$
if $\kappa_l <0$.
\end{proposition}

%%NEW a bit below
In order to describe the asymptotics of solutions as $p \to p_l$,
we apply the Lyapunov-Schmidt method (\cite[Ch.~8]{KrasZ})  to
equation (\ref{Inteq}) with the operator ${\bf A}$  being
differentiable at $0$. Since under the assumptions of Proposition
\ref{PBif2} the kernel $E_0={\rm ker\,}{\bf A}'(0)={\rm
Span\,}\{\psi_l\}$ is one-dimensional, denoting by $E_1$ the
complementary (orthogonal to $\psi_l^*$) invariant subspace, we
set $V= V_0+V_1$, where $V_0= \e_l \psi_l \in E_0$ and
 $$
  \mbox{$
 V_1=
\sum_{k \not = l} \e_k \psi_k \in E_1.
 $}
 $$
  Let $P_0$ and $P_1$,
$P_0+P_1=I$, be projections onto $E_0$ and $E_1$ respectively.
 Projecting (\ref{Inteq}) with $n=\infty$ onto $E_0$ yields
 \beq
 \label{eqInt}
  \mbox{$
\g_l \e_l =  \langle \hat \BB^{-1} (|V|^{p-1}V), \psi_l^* \rangle,
 \quad \g_l =1-\frac{ 1+c_1}{1+\frac l{2m}}= -\mbox{$\frac {(N+l)s}{(1+\a)(p-1)(2m+l)}$},
  $}
 \eeq
 where $s=p_l-p$.
By the general bifurcation theory (see e.g. \cite[p.~355]{KrasZ},
 \cite[p.~383]{Deim}, and branching approaches in \cite{VainbergTr};
  note that operator ${\bf A}'(0)$ is
Fredholm of index zero), the equation for $V_1$ can be solved and
this gives  $V_1 = o(\e_l)$ as $\e_l \to 0$, so that $\e_l$ is
calculated from the Lyapunov bifurcation equation (\ref{eqInt}) as
follows
 $$
 \g_l \e_l = |\e_l|^{p-1}\e_l \langle \hat \BB^{-1} |\psi_l|^{p-1}\psi_l, \psi_l^* \rangle
 + o(|\e_l|^p) \,\, \Longrightarrow \,\,|\e_l|^{p-1} = \hat c_l [(p_l - p)+
 o(1)],
 $$
 where $\hat c_l =
\mbox{$\frac {(N+l)^2}{4m^2(1+\a) \kappa_l}$}.
 $
 %%NEW
 Here, we have performed calculations as follows:
  $$
  \langle \hat \BB^{-1} \psi_l^p, \psi_l^* \rangle =
 \langle \psi_l^p,  (\hat \BB^*)^{-1} \psi_l^* \rangle =
 -\kappa_l/(1+\mbox{$ \frac l{2m}$}).
 $$

It is natural to require $\kappa_l >0$, though this is not always
true for  $\a>0$; see a closed
 curve of solutions in Figure \ref{FF2} (for $\a=0$ the branches
 were detected to be always monotone, \cite{GW2}).
 In view of the
orthonormality property (\ref{Ortog}), for $p=1$ we have
$\kappa_l=1$, so that by continuity we can guarantee that
 \beq
 \label{kap+}
\kappa_l > 0 \quad \mbox{at least for all} \,\,\, p \approx 1^+.
 \eeq
 %%NEW
Thus, we obtain a countable sequence of bifurcation points
(\ref{pcrit}) satisfying $p_l \to 1^+$ as $l \to \infty$, with
typical pitch-fork bifurcation branches appearing in a left-hand
neighbourhood, for  $p < p_l$. The  behaviour of solutions in
$H^{2m}_\rho$ and uniformly takes the form
 \beq
  \label{Vexpn}
 V_l(y) = \pm \left[ \hat c_l (p_l - p) \right]^{1/(p-1)}
( \psi_l(y) + o(1)) \quad \mbox{as} \,\,\, p \to p_l^-.
 \eeq

We now prove the main result concerning ``local" existence and
stability of the VSS solution with the  similarity profile
$V_0(y)$ corresponding to the first bifurcation point, $p = p_0$.
If $\kappa_0 > 0$, as expected, then two bifurcation branches
exist for $p <p_0$.

\begin{theorem}
\label{ThEx} For  $p \approx p_0^-$,  the problem $(\ref{1.7R})$
 %%%$(\ref{ExpC})$
 admits a solution $V_0(y) \not \equiv 0$ provided that
$\frac{2m}N$ is small enough, and then it is an asymptotically
stable stationary solution.
 %% of the rescaled equation
 %%%$(\ref{veq1})$.
\end{theorem}

\noi {\em Proof.} As we have shown,
% monotone
% decreasing (in $p$)
a continuous branch bifurcating at $p=p_0^-$ exists if
 \beq \label{K00}
  \mbox{$
 \kappa_0 = \langle |\psi_0|^{p_0-1}\psi_0, \psi_0^* \rangle \equiv
\int |F|^{\frac{2m}N} F \, > 0 \quad (\psi_0^* \equiv 1).
 $}
 \eeq
  In view of
 the
positivity dominance of the rescaled  fundamental solution  $F$,
$\int F = 1$,
 we have that (\ref{K00}) holds by continuity provided that $\frac{2m}N
\ll 1$. Therefore, in this case there exists a solution
(\ref{Vexpn}) with $l = 0$ satisfying for small $s = p_0 - p > 0$
uniformly
 \beq
 \label{V01}
  V_0(y) = (\hat{c}_0 s)^{\frac{1}{p-1}}[F(y)+o(1)],
   \quad \hat{c}_0 = \mbox{$\frac {N^2}{4m^2 \kappa_0}$}.
 \eeq
We now estimate  the spectrum of the corresponding linearized operator
 %%% of
 %%equation (\ref{veq1})
 \beq
\label{BBC}
    {\bf D}_0 = \BB_1 - p |V_0|^{p-1} I.
 \eeq
 %%%NEW
Some of the eigenvalues of the operator (\ref{BBC}) follow from
the original PDE (\ref{1.1R}). For instance, the stable eigenspace
with
 $
 \hat \l= -1$, $ \hat \psi =  \frac 1{p-1}V_0 + \frac 1{2m}
y \cdot  \nabla V_0 \in L^2_\rho,
 $
 follows from the time-translational invariance of the PDE. For
 $N=1$, translations in $x$ yield another pair
  $
  \hat \l = - \frac 1{2m}$, $\hat \psi =
  V_{0y} \in L^2_\rho.
  $
 For $N>1$, in the non-radial setting, this $\hat p$ has multiplicity $N$ with eigenfunctions
  $V_{0 y_i}$. These are not the first pair with
%the maximal ${\rm Re} \, \hat p$.
the maximal real part.

Bearing in mind that the spectrum of the unperturbed operator
$\BB$ is real, (\ref{spec1}), and has the unique, non-hyperbolic
eigenvalue $\lambda_0 = 0$, we use (\ref{V01}) to obtain
 \beq
 \label{BBC1}
 {\bf D}_0 =  \BB + s(1+o(1)) {\bf C},
\eeq
  where, as it follows from  (\ref{K00}) and (\ref{V01}) at $p=p_0$,
 the perturbation has the form
   \beq
\label{CCC}
   {\bf C} = \mbox{$\frac{N^2}{4 m^2(1+\a)} \big(1 - \frac{p_0}{\kappa_0} |F|^{\frac {2m(1+\a)}N} \big)I.$}
 \eeq
Therefore, we consider the spectrum of the perturbed operator
  \beq
  \label{BCCB}
 \tilde {\bf D}_0= \BB + s {\bf C}.
 \eeq
 Since $(\BB-I)^{-1}{\bf C}$ is bounded,
 $$
 (\tilde {\bf D}_0-I)^{-1} = [I +s(\BB-I)^{-1}{\bf C}]^{-1}(\BB -I)^{-1}
 $$
  is compact for small $|s|$ as the product of a compact and bounded operators.
  Hence,
 $\tilde {\bf D}_0$ also has only a discrete spectrum.
By the classical perturbation theory of linear operators (see e.g.
\cite{GohKr}), the eigenvalues and eigenvectors of $\tilde {\bf
D}_0$ can be constructed as a perturbation of the discrete
spectrum $\sigma(\BB)$ consisting of eigenvalues of finite
multiplicity. We are interested in the perturbation of the first
simple eigenvalue $p_0=0$. Setting
\[
\tilde{\lambda}_{0} =  s \mu_{0} + o(s), \; \; \tilde{\psi}_0 =
\psi_0 + s \varphi_0 + o(s) \quad \mbox{ as }\,\, s \to 0
\]
and substituting these expansions in the eigenvalue equation $
\tilde {\bf D}_0 \tilde{\psi}_0 = \tilde{\lambda}_{0}
\tilde{\psi}_0
$
yields
\beq
\label{BBtilde}
 \BB \varphi_0 =
(-{\bf C} + \mu_0 I) \psi_0.
 \eeq
We then obtain the solvability (orthogonality) condition
\[
\langle (-{\bf C} + \mu_0 I) \psi_0,\psi_0^* \rangle = 0 \,\,
\Longrightarrow \,\, \mu_0 = \langle {\bf C}F, 1 \rangle.
\]
Using  (\ref{CCC}) yields
 $$
 \mbox{$
  \mu_0 = -  \frac N{2m} < 0.
   $}
   $$
Therefore, ${\rm Re}\, \tilde{\lambda}_0 < - \frac {N s}{4m} <0$
for all $p \approx p_0^-$. Since, with these properties of the
spectrum, the perturbation (\ref{BBC}) of $\BB$ remains a
sectorial operator with $\s(\tilde {\bf D}_0) \subset \{{\rm Re}
\, \l \le - \frac  {N s}{4m}\}$ and $\|{\mathrm e}^{\tilde \BB_1
\t}\|_{\mathcal L} \le C {\mathrm e}^{-N(p_0-p) \t/4m}$ in the
norm of ${\mathcal L}(H^{2m}_\rho,H^{2m}_\rho)$ \cite{Fr},
$V_0(y)$ is exponentially stable in $H^{2m}_\rho $. $\qed$

 \vspace{.2cm}

  We expect that
condition (\ref{K00}) remains valid for any $m$ and $N$ so that
$V_0(y)$ is stable without the restriction $2m \ll N$.  We have a
numerical support for this, but, as yet, no rigorous proof exists.
%%%NEW
Possibly, to check conditions such as (\ref{K00}) we must
currently rely on numerical evidence and then, as often happens in
spectral theory and applications, Theorem \ref{ThEx} can be
established with a hybrid analytic-computational proof.
 We also expect that
the whole branch bifurcating from $p=p_0$ remains stable for all
$p \in (1,p_0)$, though the proof would require to establish that
the discrete spectrum $\s({\bf D}_0)$ never touches the imaginary
axis. In particular, this  open problem means that a new
(nonlinear) saddle-node bifurcation never occurs on this
$p_0$-branch, i.e., it does not have turning points.
%%NEWVAG
For the variational problem with $m=1$, this is valid  \cite{SW}
as well as for ordinary differential higher-order equations with
self-adjoint positive operators of special structure of
quasi-derivatives \cite{Rynn2, Rynn1}; see also  properties of
bifurcation branches in  Berger \cite[p.~380]{Berger}.

Further, it is easy to see that  the other bifurcation branches
are {\em unstable}. Taking  any $l \ge 1$, instead of (\ref{BBC1})
we now have
\[
 {\bf D}_l = \BB_1 - p |V_l|^{p-1} I \equiv \BB + \big[c_1 -
  s p_l \hat c_l (|\psi_l|^{p-1} +o(1))\big]I, \quad s= p_l - p.
\]
From the definition of $\BB_1$, (\ref{BB*}), $c_1 > 0$ for all $p
\approx p_l$, thus  $V_l$ for $l \ge 1$ is unstable.

%%%%%%%%%%%%%%%%%%%%%%%%%%%%%%%%%%%%%%%
\subsection{Remark on global bifurcation diagrams}

Global bifurcation results concerning continuous branches of
solutions  were already given in Krasnosel'skii
\cite[p.~196]{Kras} (the first Russian edition was published in
1956). Concerning further results and extensions, see references
in \cite[Ch.~10]{Deim} (especially, see \cite[p.~401]{Deim} for
typical global continuation of bifurcation branches), and also
\cite[\S~56.4]{KrasZ}. These approaches deal with integral
equations with compact operators such as (\ref{Inteq}).

 In the present non-variational problem, the main open
problem of concern is to establish under which conditions the
$p$-branch originated at some $p=p_l$ can be extended up to
$p=1^-$, so that cannot end up at another bifurcation point $p_k <
p_l$. Actually, this can happen for $\a>0$; see Figure \ref{FF2},
where the closed branch is originated at $p=p_0$ and ends up at
$p=p_2$.

%%%%%%%%%%%%%%%%%%%%%%%%%%%%%%%%%%%%%%%%%%%%%%%%%%%%%%%%%
%%%%%%%%%%%%%%%%%%%%%%%%%%%%%%%%%%%%%%%%%%%%%%%%%%%%%%
\section{%%%On spatial structure of various VSSs profiles:
$p$- and $\a$-bifurcations and branches:
%% analytic and
%%%%numerical evidence of
 VSS multiplicity}
 \label{SectAl}

We consider the fourth-order equation (\ref{1.7R}) in 1D with
$m=2$, i.e., the ODE problem
 \beq
 \label{4s}
 %% \mbox{$
  \left\{
  \begin{matrix}
 -V^{(4)} + \frac 14 \, y V' + \b V -|V|^{p-1}V=0 \quad \mbox{in}
 \quad \re \quad \big(\b= \frac{1+\a}{p-1}\big),\smallskip\smallskip\\
 V(y) \,\,\, \mbox{has exponential decay as $y \to \infty$}.
 \qquad\qquad\qquad\quad
%% $}
 \end{matrix}
  \right.
  \eeq
Note that the condition of the exponential decay is crucial for
VSS setting, since the ODE in (\ref{4s}) admits continuous
families of solutions with algebraic power decay such as
  \beq
  \label{bb1}
 V(y) = C_\pm |y|^{-\frac{2m}{p-1}}(1+o(1)) \quad \mbox{as} \quad y \to
 \pm \infty,
  \eeq
  where $C_\pm \in \re$ are, in general, arbitrary constants; see
  an example in
  \cite[\S~8]{GW2} for the case $\a=0$. In numerical
  applications, the exponentially decaying solutions of
  (\ref{4s}) are always clearly oscillatory for $y \gg 1$ that
  strongly differ them from those satisfying  (\ref{bb1}).

 For even
profiles $V_0$, $V_2$,... satisfying (\ref{4s}) that we pay more
attention to, we impose the symmetry conditions at the origin
 \beq
 \label{Symm1}
 V'(0)=V'''(0)=0.
 \eeq
For the odd profiles $V_1$, $V_3$,...\,, we pose the anti-symmetry
conditions
 \beq
 \label{Asymm1}
 V(0)=V''(0)=0.
  \eeq

%%%%%%%%%%%%%%%%%%%%%%%%%%%%%%%%%%%%
\subsection{$p$-branches for fixed $\alpha \ge 0$}

\begin{figure}
%\vskip -.3cm
\centering \subfigure[$p_0$-branch]{
\includegraphics[scale=0.5]{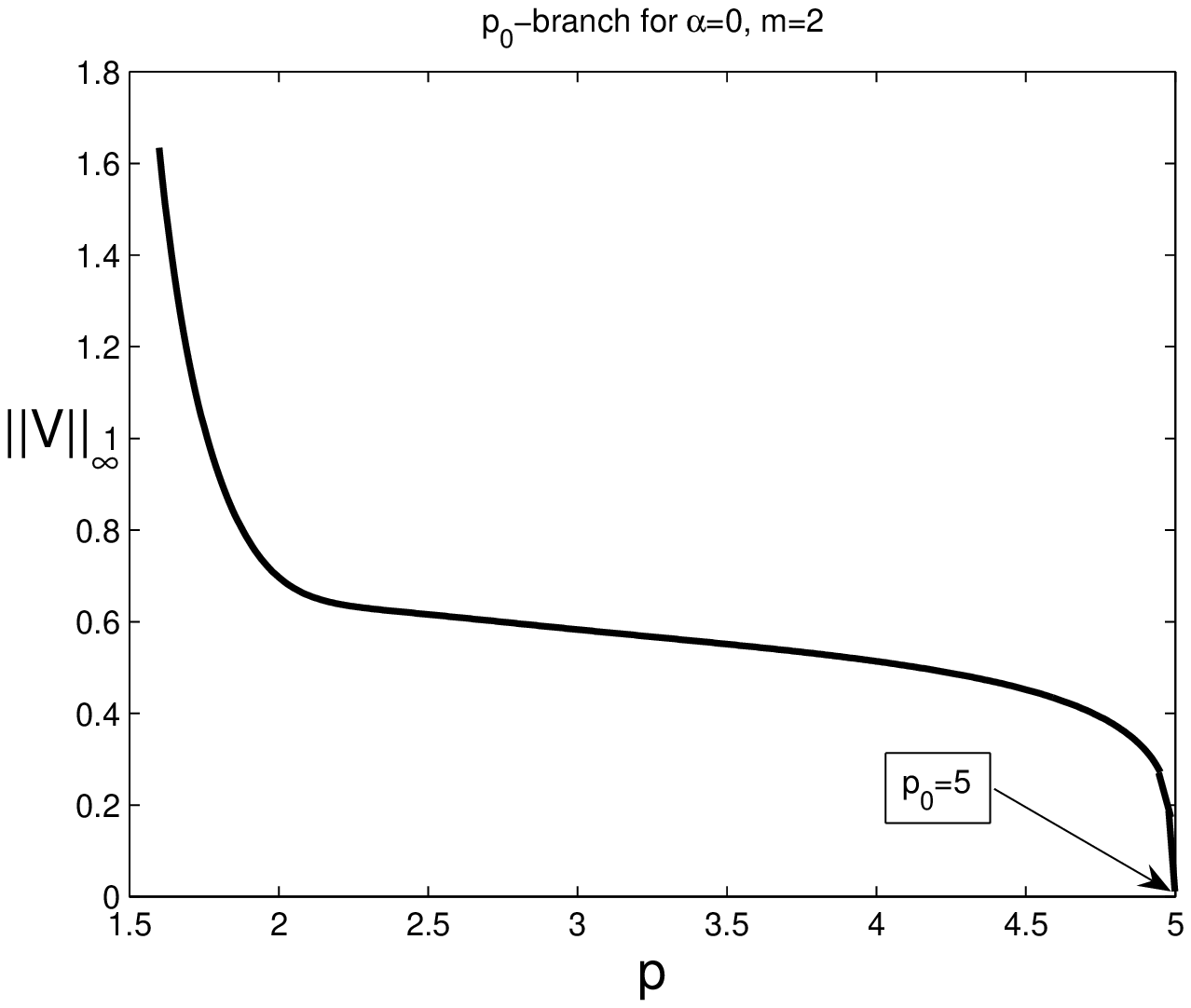}                %%{m2a0.eps} %Old: no N
} \subfigure[$V_0$ profiles]{
\includegraphics[scale=0.5]{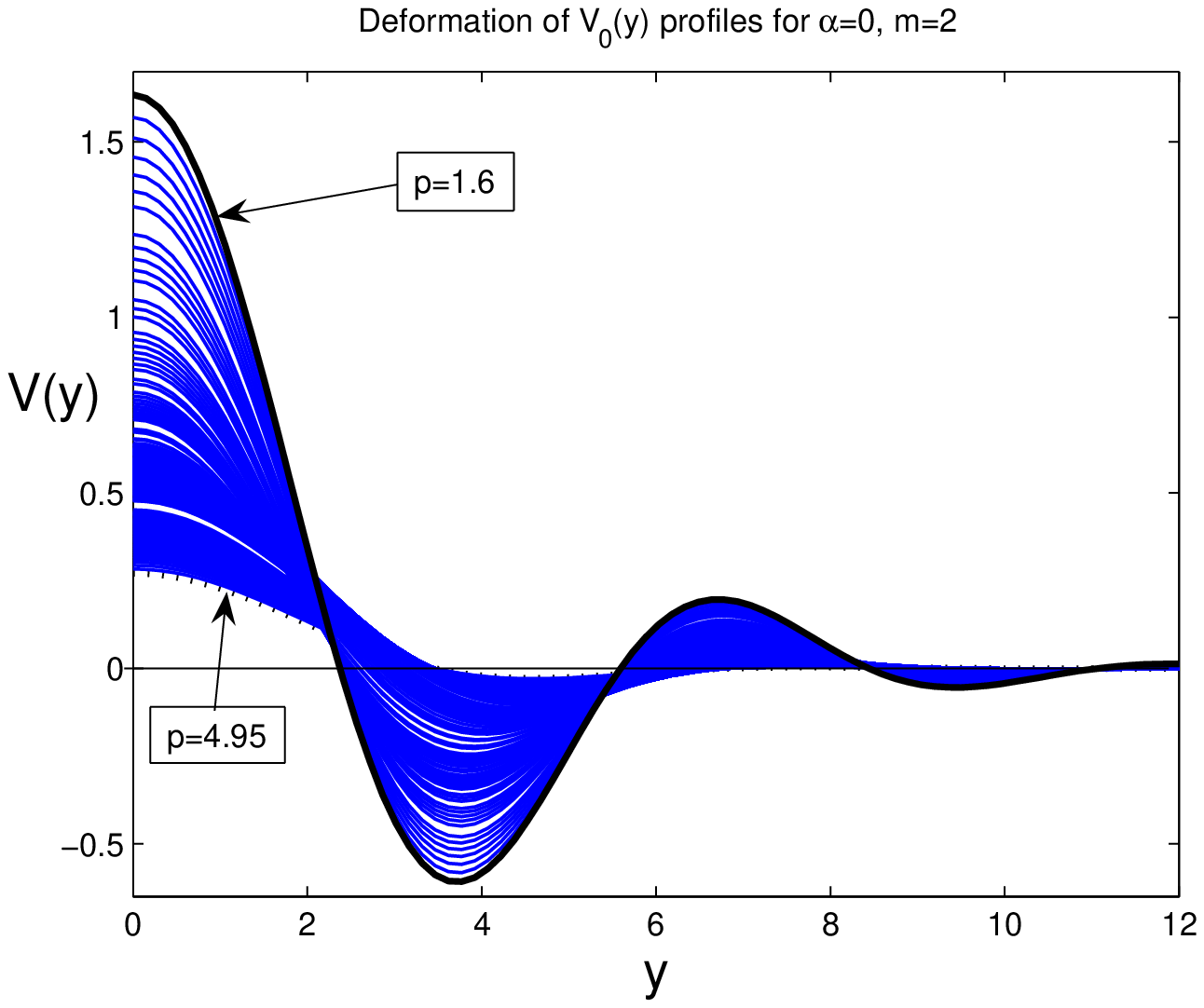}                  %%{m2a0Pr.eps}
}
 \vskip -.4cm
\caption{\rm\small The first $p_0$-branch of solutions of
(\ref{4s})  for  $m=2$ and $\a=0$; the branch (a) and deformation
of profiles $V_0(y)$ with $p$, (b).}
  \vskip -.1cm
 %of the
 %ODE (\ref{pp.5}) for $m=-0.75$ (a) and $m=-0.9$ (b).}
% for
 %$n=1$ near non-oscillatory (a) and oscillatory (b) interfaces.}
 \label{FF1}
\end{figure}

 As we have mentioned, the autonomous case $\a=0$ is
rather well-understood, \cite{GW2}. For example, in Figure
\ref{FF1}, we present the first strictly monotone $p_0$-branch of
similarity profiles $V_0(y)$ for $\a=0$. This branch blows up as
$p \to 1^-$ according to the asymptotics in \cite[p.~1091]{GW2}.

%%%%%%%%%%%%%%%%%%%%%%%%%%%%%%%%%%%%%%%%%%%%%%%%%%
%%%%%%%%%%%%%%%%%%%%%%%%%%%%%%%%%%%%%%%%%%%%%%%%%%%%%%%%%%%%%%%%%%%%%%%%%

\begin{figure}
%\vskip -.3cm
\centering \subfigure[$p_0$-branch]{
\includegraphics[scale=0.5]{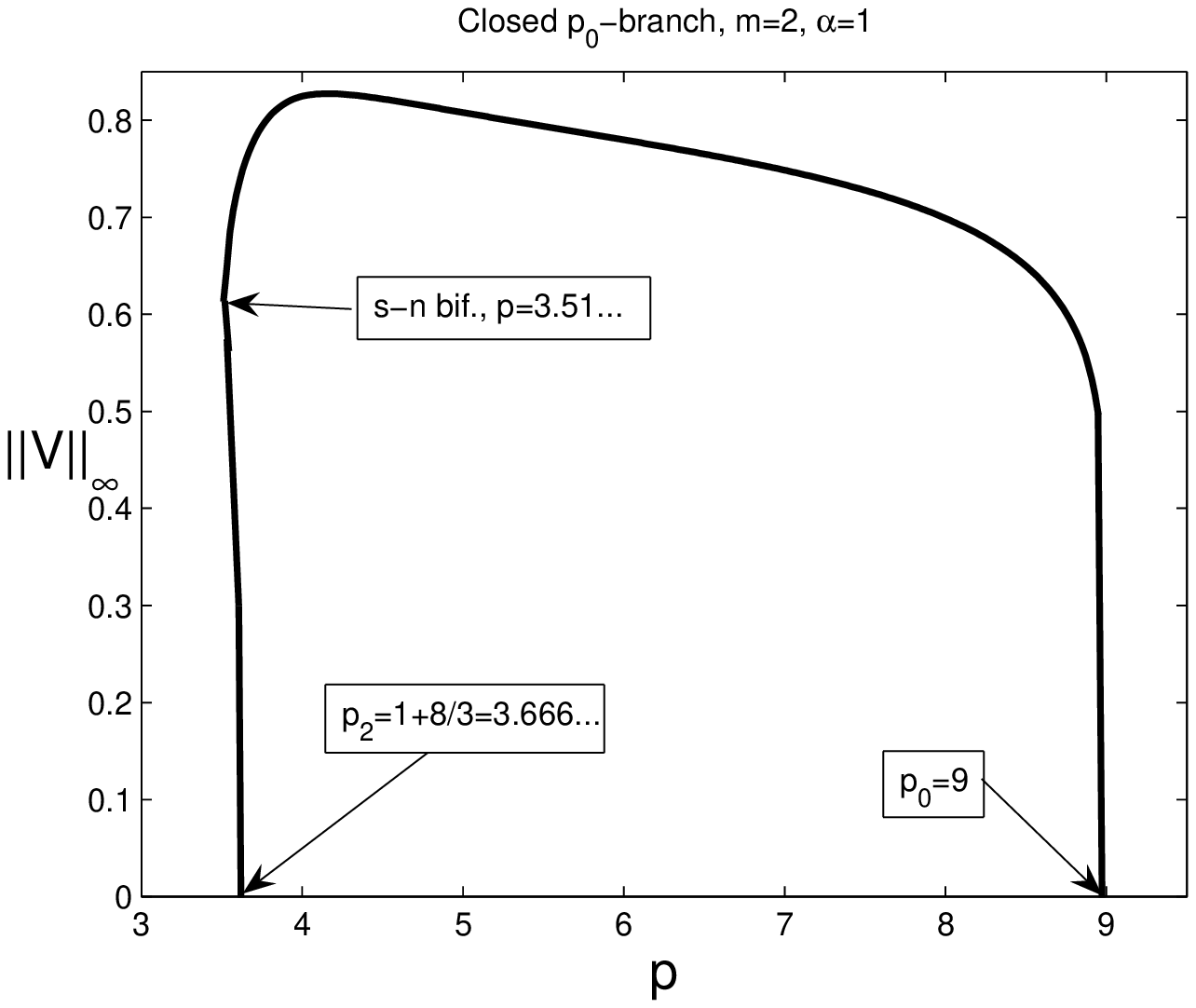} %Old: no N
} \subfigure[$V_0$ profiles]{
\includegraphics[scale=0.5]{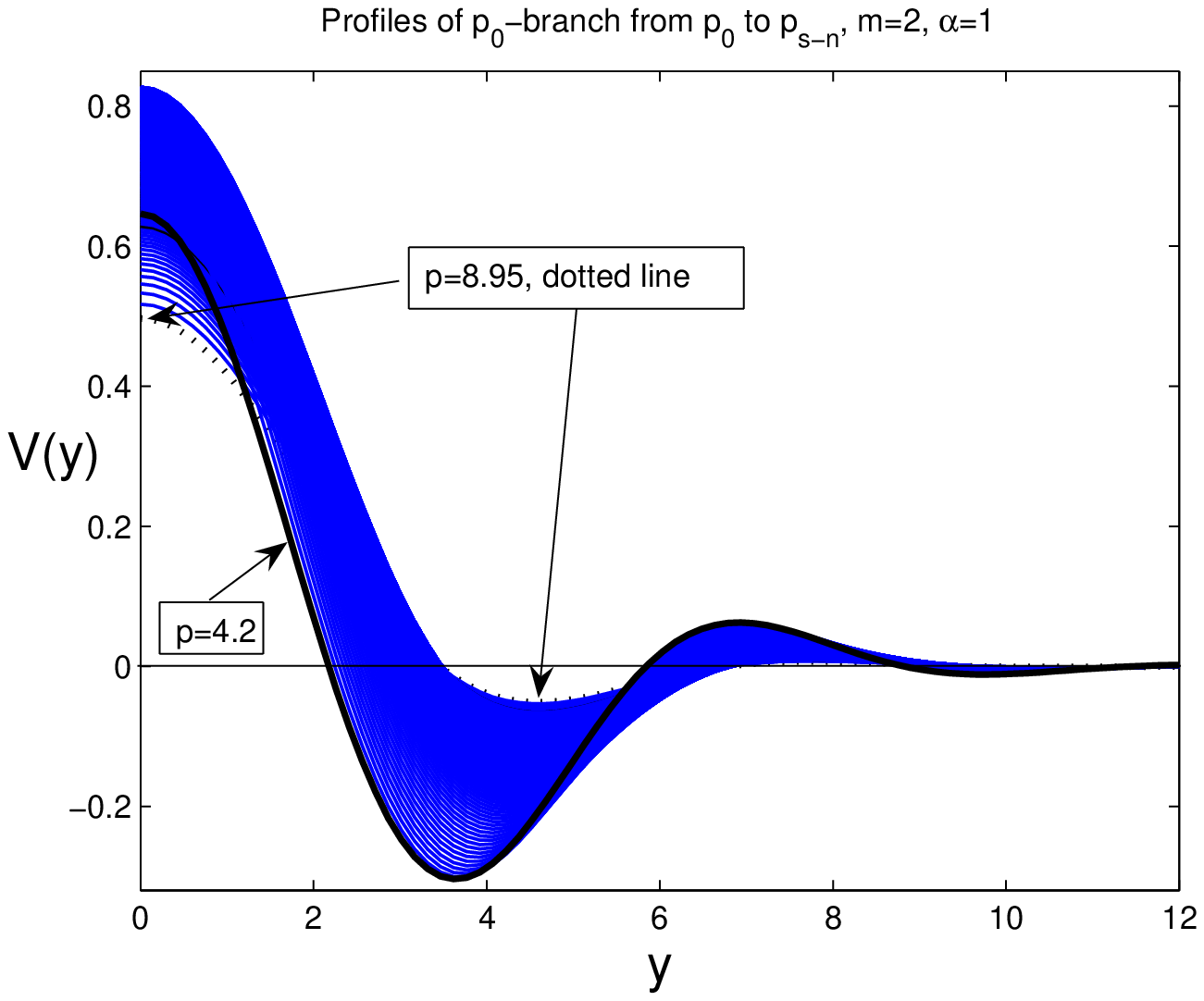}
}
 \vskip -.4cm
\caption{\rm\small The closed  $p_0$-branch of solutions of
(\ref{4s}) for  $m=2$, and $\a=1$; the branch (a) and deformation
of profiles $V_0(y)$, (b).}
  \vskip -.1cm
 %of the
 %ODE (\ref{pp.5}) for $m=-0.75$ (a) and $m=-0.9$ (b).}
% for
 %$n=1$ near non-oscillatory (a) and oscillatory (b) interfaces.}
 \label{FF2}
\end{figure}

\begin{figure}
%\vskip -.3cm
\centering \subfigure[$p_0$-branch]{
\includegraphics[scale=0.5]{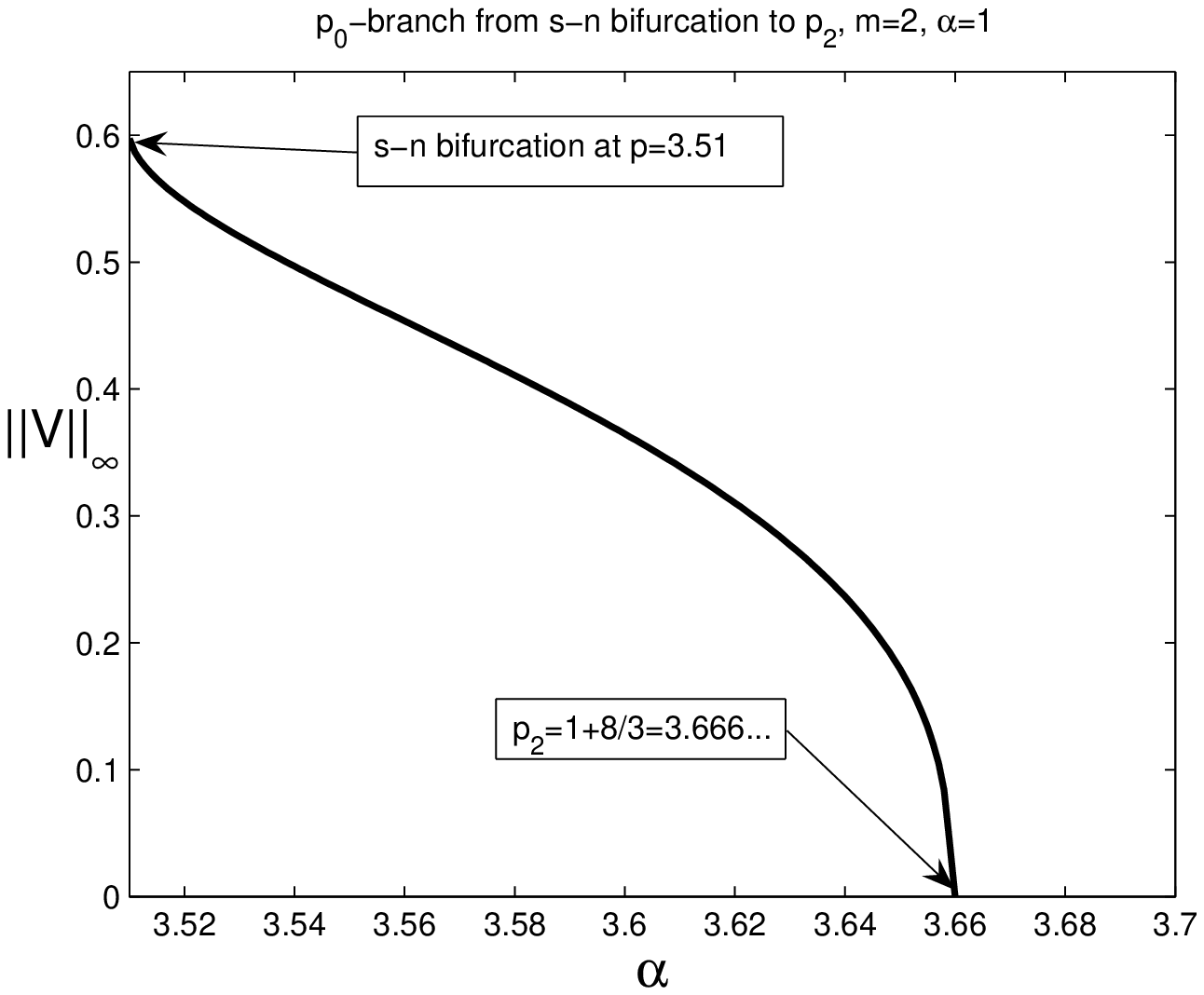} %Old: no N
} \subfigure[$V_0$ profiles]{
\includegraphics[scale=0.5]{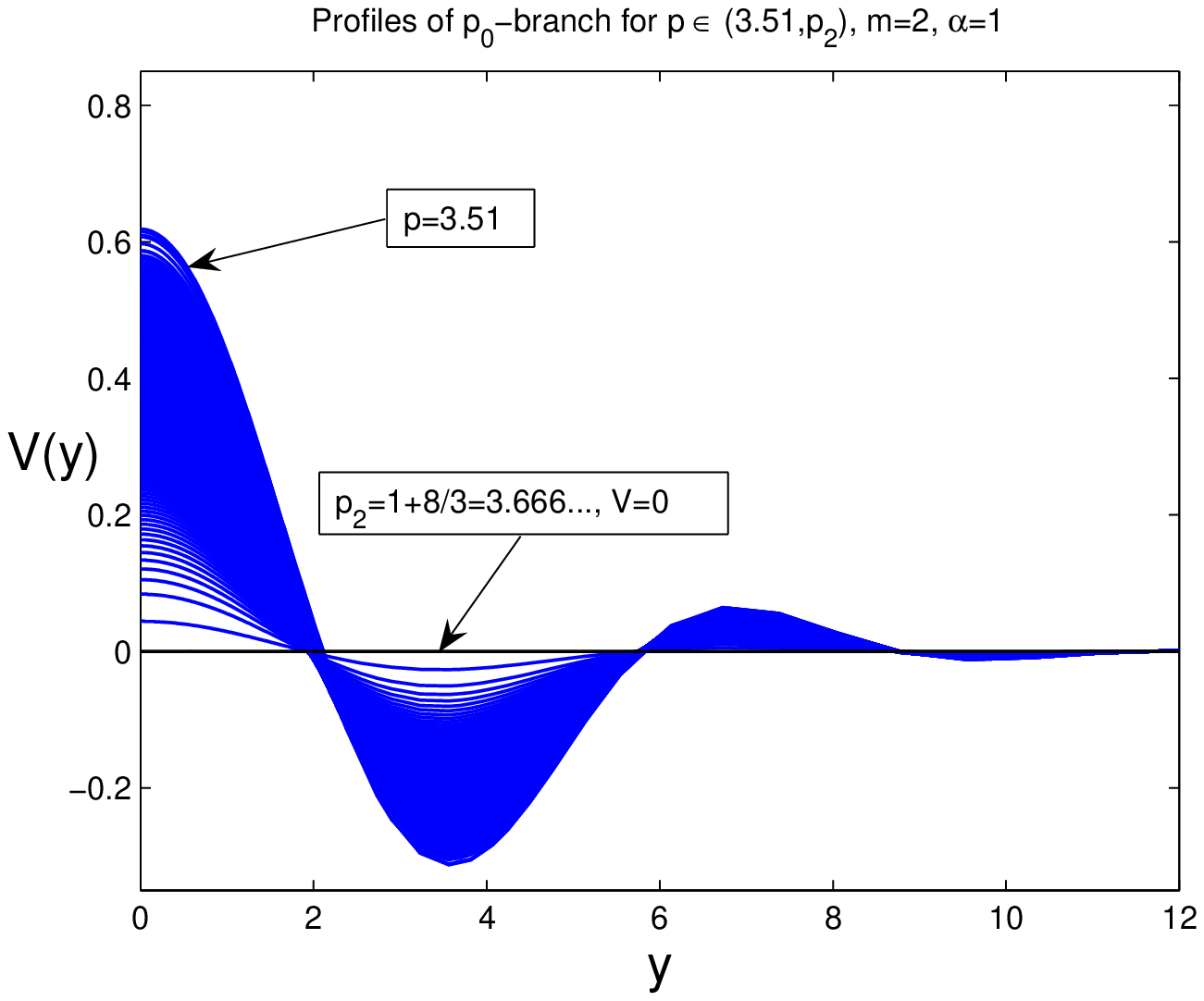}
}
 \vskip -.4cm
\caption{\rm\small The   $p_0$-branch of solutions of (\ref{4s})
between two bifurcations;
 $p \in [3.51,3.66]$, $m=2$, $\a=1$; the
branch (a) and deformation of profiles $V_0(y)$, (b).}
  \vskip -.1cm
 %of the
 %ODE (\ref{pp.5}) for $m=-0.75$ (a) and $m=-0.9$ (b).}
% for
 %$n=1$ near non-oscillatory (a) and oscillatory (b) interfaces.}
 \label{FF2N}
\end{figure}

 For $\a>0$, the $p$-bifurcation
 branches can essentially change their topology and can be  closed
 curves, as Figure \ref{FF2} demonstrates for $\a=1$.
 Here, the first $p_0$-branch is initiated at $p_0=9$ and has the
 {\em saddle-node} (turning) point at
  $$
  p_{\rm s-n}=3.51... \, .
  $$
Eventually, the branch ends up at the third critical point
(\ref{bif1}), i.e.,
%%% for the given exponents,
  at
 $$
 \mbox{$
 p_2= 1 + \frac 83=3.666... \, .
  $}
  $$
  For convenience, in Figure \ref{FF2N}, we give the enlarged
  structure of the almost vertical part of the branch in
  Figure \ref{FF2}(a) and VSS profiles on it for $p \in[p_{\rm
  s-n},p_2]$.

%%%%%%%%%%%%%%%%%%%%%%%%%%%%%%%%%%%%%%%%%%%%%%%%%%
%%%%%%%%%%%%%%%%%%%%%%%%%%%%%%%%%%%%%%%%%%%%%%%%%%%%%%%%%%%%%%%%%%%%%%%%%

%%%%%%%%%%%%%%%%%%%%%%%%%%%%%%%%%%%%%%%%%%%%%%%%%%
%%%

%%%%%%%%%%%%%%%%%%%%%%%%%%%%%%%%%%%%%%%%%%%%
 \subsection{On $\a$-branches for fixed $p$}

\begin{figure}
%\vskip -.3cm
\centering \subfigure[$\a_0$-branch]{
\includegraphics[scale=0.5]{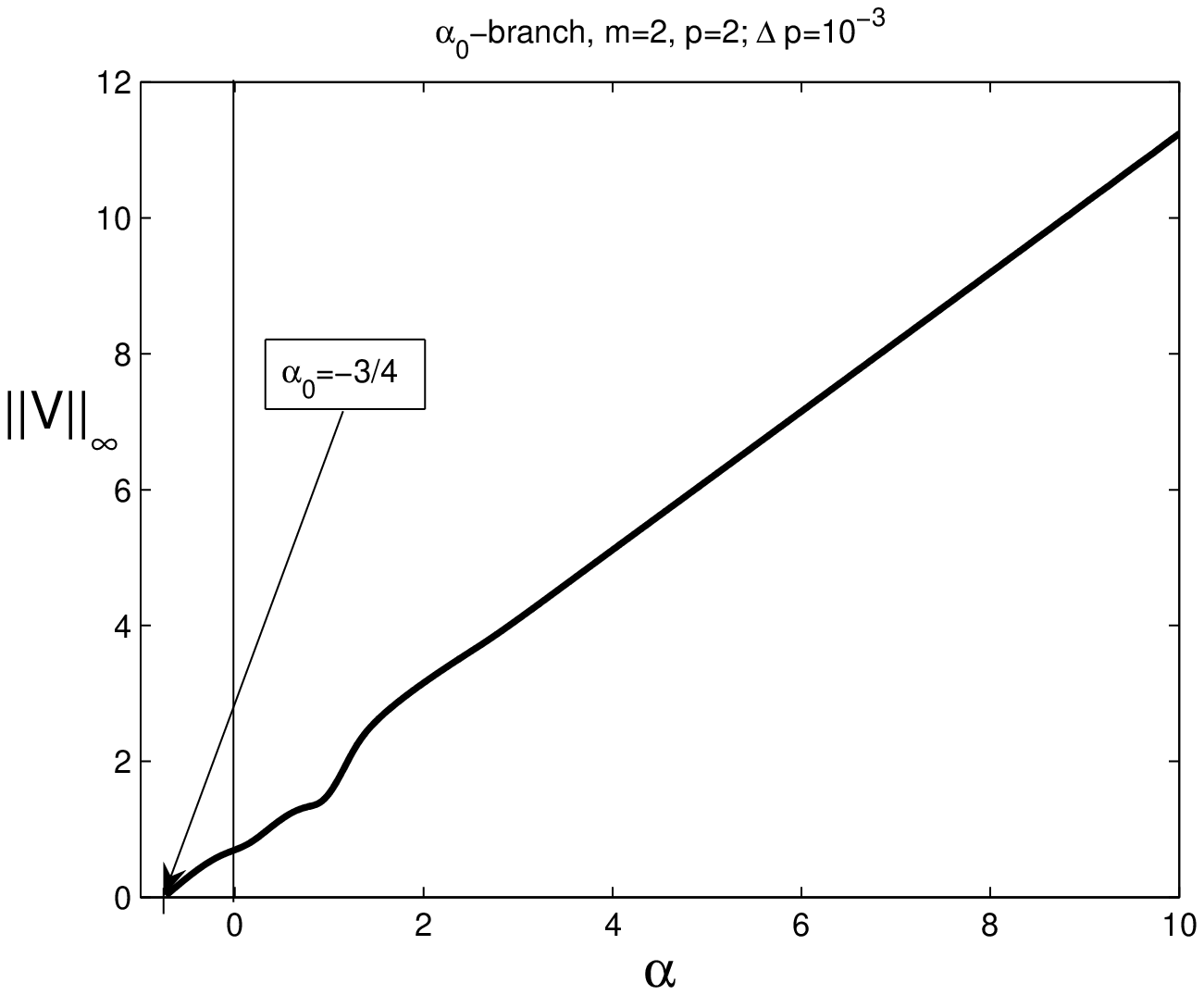} %Old: no N
} \subfigure[$V_0$ profiles]{
\includegraphics[scale=0.5]{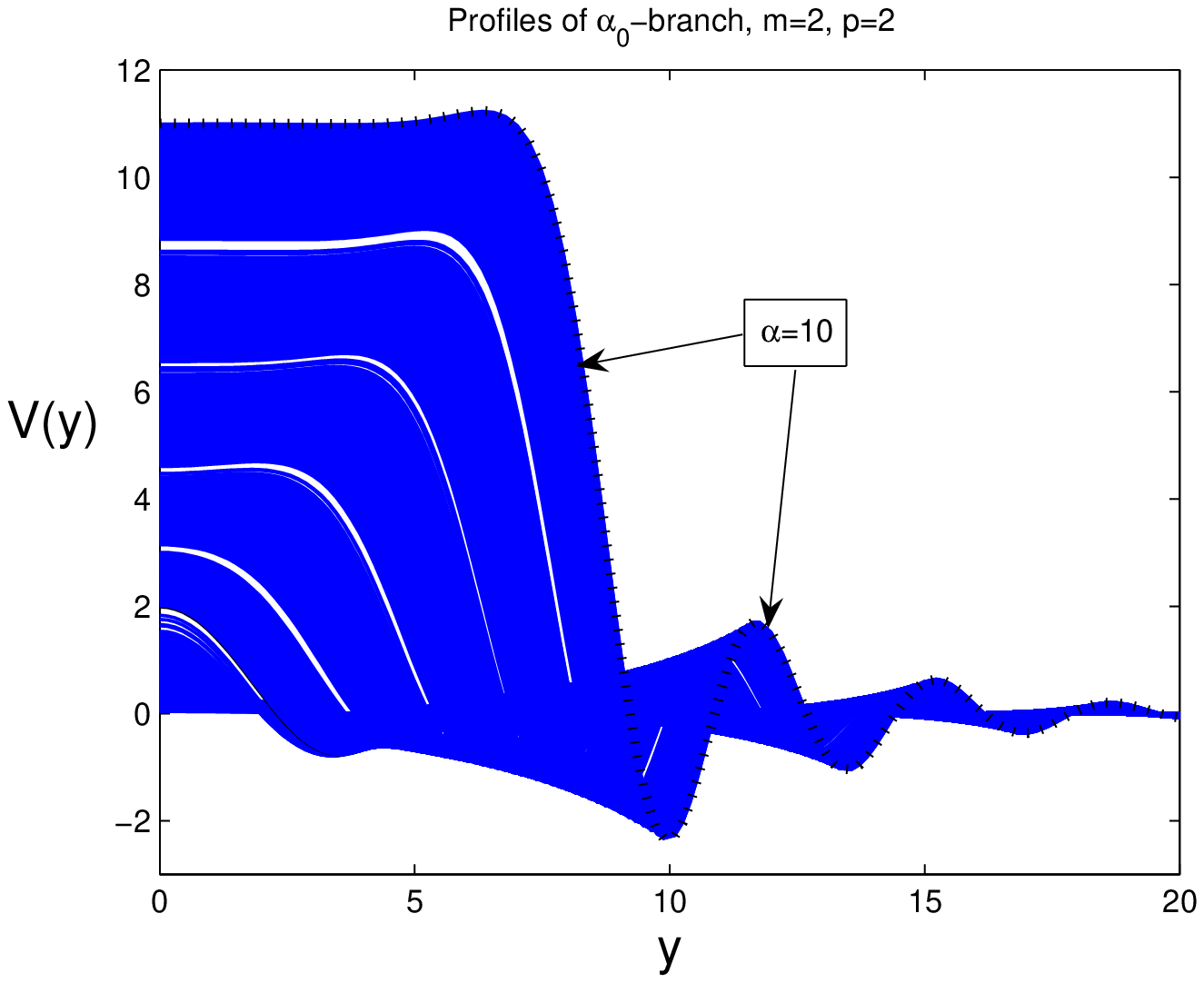}
}
 \vskip -.4cm
\caption{\rm\small The first  $\alpha_0$-branch of solutions of
(\ref{4s}) for  $m=2$, and $p=2$; the branch (a) and deformation
of profiles $V_0(y)$, (b).}
  \vskip -.1cm
 %of the
 %ODE (\ref{pp.5}) for $m=-0.75$ (a) and $m=-0.9$ (b).}
% for
 %$n=1$ near non-oscillatory (a) and oscillatory (b) interfaces.}
 \label{FF3}
\end{figure}

 A typical example is presented in Figure \ref{FF3}, where we show the first
 $\alpha$-branch for $N=1$, $p=2$, that, according to (\ref{al11}), is originated at
 the bifurcation point
  $$
   \mbox{$
  \a_0=-\frac 34 \quad (l=0).
   $}
  $$
We observe from (b) that the profiles get wider as $\a$ increases.
This is a natural phenomenon that will play an important role in
what follows. It should be noted that  Figure (b) does not explain
jump-discontinuities that are not only possible but often happen
in such numerics (in view of huge density of various branches of
solutions), and were observed even with the enhanced Tolerances
$10^{-4}$ and the step $\Delta \a=10^{-3}$. In fact, (a) shows for
$\a \approx 1$ a certain non-smoothness of the branch that we
cannot explain. To be honest, we cannot guarantee that here this
is not associated with a discontinuous jump to other neghbouring
branches. The $\a$-branching deserves both extra numerical and
analytical study.

%%In Figure \ref{FF3N}, we present the enhanced variant of these
%%numerics, where (b) clearly shows discontinuities that occur in
%%view of closedness of various $\a$-branches of neighbouring VSS
%%profiles.

%%%%%%%%%%%%%%%%%%%%%%%%%%%%%%%%%%%%%%%%%%%%%%%%%%
\subsection{Further bifurcations and discussion}

%%Therefore, each
%%function from
We denote  the basic set of VSS profiles by
 %%denoted by (sometimes we omit the
%%upper subscript)
 \beq
 \label{BS1}
\mbox{basic spectrum:} \quad  V_0, \,\,  V_1, \,\, V_2, \, ..., \,
V_k, \, ... \, .
 %% \quad  V_0^{(0)}, \,\,  V_1^{(0)}, \,\,
%%V_2^{(0)}, \, ..., \,  V_k^{(0)}, \, ... \, ,
  \eeq
%% may generate other profiles (branches) at saddle-node
%% bifurcations.
%% Note that
  Each  $V_k(y)$ satisfies the {\em approximate Sturmian
 property}, i.e.,
 has precisely $k$ {\em dominant} extrema. For $m=1$, this is true sharply, by the Maximum Principle.
 For $m \ge 2$, by ``dominant"  we mean
 those minima or maxima that are different from smaller ones in the
 oscillatory tail of the solutions.
  This approximate Sturm's zero (extremum) property for $N=1$
 is
  associated with the fact that in the ODE (\ref{1.7R}) with $m
  \ge 2$, the leading $2m$th-order operator is a multiplication
  of $m$ positive second-order ones,
   $$
   -(-1)^{m} D_y^{2m} \equiv -(-D_y^2)^m,
   $$
   for which Sturm's property is true concerning its
   eigenfunctions; see general conclusions in Ellias
   \cite{Ellias},  also applications to bifurcation diagrams of interest here
    in \cite{Rynn2, Rynn1}. It turns out that lower-order perturbations in
   (\ref{1.7R}) create only sufficiently small extrema and
   infinitely many zeros in the exponential tails of solutions.

\smallskip

\begin{figure}
%  \vskip -.3cm
\centering
\includegraphics[scale=0.65]{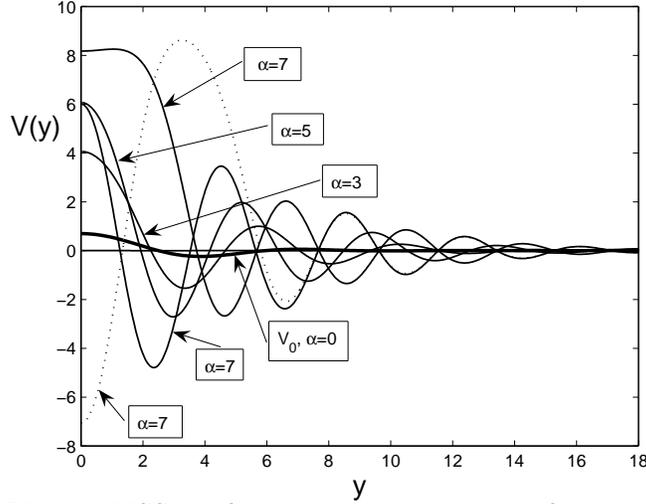}  %%%%%%%{4het5.eps}   %%%%%{4het.eps}  Old
\vskip -.5cm \caption{\small Various VSS profiles associated with
the first basic pattern $V_0(y)$ for $p=2$, $N=1$,
 and different
values of $\a \in [0,7]$.}
   \vskip -.3cm
 \label{F00}
\end{figure}

For $\a>0$, the basic profiles are essentially deformed. For
instance,
  a few of such profiles (already available in Figure \ref{FF3}(b)) are
 shown in Figure \ref{F00} for various $\a  \ge 0$.

 %% where we also present another type of
 %%VSS profiles such as $V_0^{(1)}(y)$.
 %% More precisely, we claim
 %%that
 %%  for sufficiently large $\a>0$,
 %%  there exist other VSS
%%sets  denoted by $\{V_k^{(j)}, \, k,j \ge 0\}$ which are obtained
%%from functions similar to those from the basic spectrum
%%(\ref{BS1}) by a cascade of saddle-node bifurcations at some
%%values $\a= \a_k^j$ beforehand.
%% We do not know precisely the cardinal number
%%of each set but
%%We will present some evidence that the number of such branches is
%%infinite and countable. Note that the operators in our ODE are
%%neither symmetric in their linear part, nor potential. Fixing
%%those saddle-node bifurcations analytically is a difficult open
%%problem. On the other hand, these can be detected them by hard
%%numerical evidence.

 A principal new feature of the
problem with the exponent $\a$ is that there may occur other
saddle-node bifurcations in the ODE in (\ref{4s}). The origin of
extra $\a$-bifurcations can be seen as follows. The leading
operator in (\ref{1.7R}),
 $$
 \tilde {\bf B} V = -(-\D)^m V- |V|^{p-1}V
  $$
is coercive and negative in the topology of the Sobolev space
$H^{2m}(B_l) \cap L^{p+1}(B_l)$ in the subcritical range
(\ref{sob}) (this is enough for $N=1$).
%%% $$
%%% \mbox{$
%% 1< p < \frac {N+2m}{N-2m} \quad (N>2m).
%%  $}
%%%  $$
Here $B_l$ is the ball of the radius $l>0$ in $\ren$. It follows
from Lemma \ref{lemspec}(i) that a similar result remains true in
the weighted space $L^2_\rho(\ren)$, and moreover, the first
linear term $\frac 1{2m}\, y \cdot \n V$ can be estimated via the
leading operator. We then observe that the crucial part is played
by the second  term
 $$
  \mbox{$
 \b V \equiv \frac {1+\a}{p-1} \, V.
  $}
 $$
The multiplier $ \frac {1+\a}{p-1}$ gets arbitrarily large
positive for $\a \gg 1$, and hence can produce a bifurcation,
where the negative and positive operators are sufficiently
balanced; see more precise estimates below.

%%%%%%%%%%%%%%%%%%%%%%%%%%%%%%%%%%%%%%%%%%%%%%%%%%%%%%%%%%%%
%%%%%%%%%%%%%%%%%%%%%%%%%%%%%%%%%%%%%%%%%%%%%%%%%%%%%%%%%%%%%%

\begin{figure}
%  \vskip -.3cm
\centering
\includegraphics[scale=0.65]{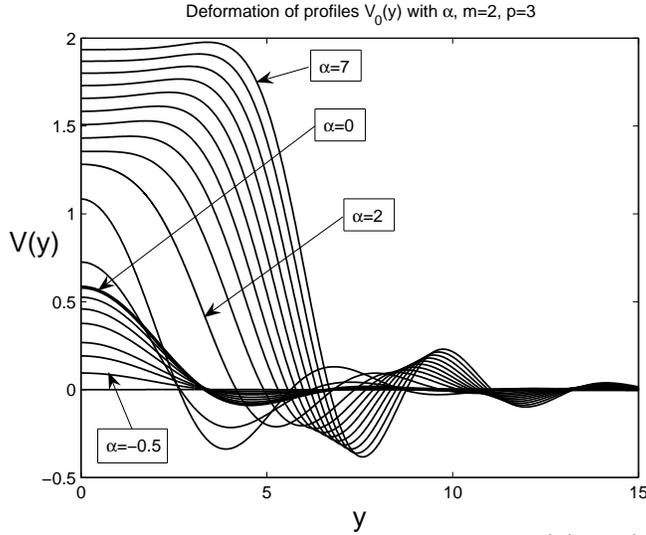}  %%%%%%%{4het5.eps}   %%%%%{4het.eps}  Old
\vskip -.5cm \caption{\small Deformation with $\alpha$ of the
first profile $V_0(y)$ of (\ref{4s}),  with $p=3$,  for $\a
 \in [-0.5, 7]$.}
   \vskip -.3cm
 \label{F1}
\end{figure}

Let us discuss THE  numerical evidence more systematically.
 In Figure \ref{F1},
we show the first VSS profiles $V_0(y)$
%%%or $V_1^{(1)}(y)$
satisfying (\ref{1.7R})
%%%, (\ref{ExpC})
 for $p=3$, $N=1$, and for a
number of different values of $\a$.
%% A saddle-node bifurcation in
%%$\a$ occurs for $\a \approx 1.5$ at which two $V_0$-branches
%%appear, $\{V_0^{(0,1)}(y)$.
 The bold line
corresponds to $V_0$ in the classic case $\a=0$ studied in
\cite{GW2}.
%%These numerical results  confirm that branches of
%%$\{V_1^{(k)}\}$ exist not only for large values of $\a \in [0,7]$
%%but also for negative values $\a \in [- \frac 12,0)$. Actually, we
%%expect that such branches exists in the whole negative
%%integrability range $n \in (-1,0)$.

%%Similar tendencies of dependence  on $\a$ are shown in Figure
%%\ref{F2} for the case $p=2$ and $N=1$ where the saddle-node
%%bifurcation occurs at about $\a \approx 1$.

%%%%%%%%%%%%%%%%%%%%%%%%%%%%%%%%%%%%%%%%%%%%%%%%%%%%%%%%%%%%

%%%%%%%%%%%%%%%%%%%%%%%%%%%%%%%%%%%%%%%%%%%%%%%%%%%%%%%%%%%%
%%\begin{figure}
%  \vskip -.3cm
%%\centering
%%\includegraphics[scale=0.65]{V18.eps}  %%%%%%%{4het5.eps}   %%%%%{4het.eps}  Old
%\includegraphics{F1.1}
%%\vskip -.5cm \caption{\small The first solution $V_0^{0)}(y)$ of
%%(\ref{1.7R}), (\ref{ExpC}), $p=2$, $N=1$ for $\alpha \approx 0$
%%splits into $V_0^{(1)}$ for larger $\a \in (3,4]$.}
%%   \vskip -.3cm
%% \label{F2}
%%\end{figure}

\begin{figure}
%  \vskip -.3cm
\centering
\includegraphics[scale=0.65]{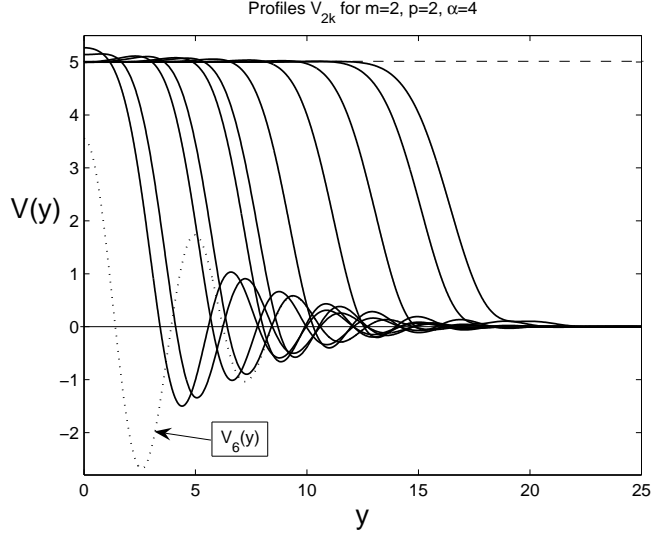}  %%%%%%%{4het5.eps}   %%%%%{4het.eps}  Old
\vskip -.5cm \caption{\small Twelve even VSS profiles $V_{2k}(y)$
%%%%%obtained from $V_0^{(0)}$
 %%%%satisfying (\ref{1.7R}), (\ref{ExpC})
for $p=2$, $N=1$, and $\a =4$.}
   \vskip -.3cm
 \label{F3}
\end{figure}

In Figure \ref{F3}, we show twelve different even VSS profiles
$V_{2k}(y)$ for $p=2$, $N=1$, and $\a=4$. Note that according to
 critical exponents (\ref{al11}), there exist precisely ten
 critical values below 4, at which such profiles can be
 originated,
  $$
   \mbox{$
  \a_0=-\frac 34, \,\, \a_2=-\frac14, \,\, \a_4=\frac14,
  \,\,...\,\,
   , \,\a_{16}=\frac {13}4, \,\, \a_{18}= \frac {15}4.
    $}
    $$
On the other hand, according to (\ref{bif1}), there exist the same
number of critical exponents $p_l$ that are above 2,
 $$
  \mbox{$
  p_0=21, \,\, p_2=\frac {23}3, \,\, ... \,\, , \,  p_{16}=
  \frac{37}{17}, \,\, p_{18}= \frac{39}{19}.
   $}
    $$
It follows that at least two profiles in Figure \ref{F3} are not
generated by standard $p$-  or $\a$-curves originated at
bifurcation points. This mystery remains an open problem and needs
extra analysis. We expect that some pairs of these profiles could
be originated at saddle-node bifurcations at some $\a_{\rm
s-n}>0$.
 It is seen that the
profiles slightly ``oscillate" about the constant equilibrium of
the ODE in  (\ref{4s}) which is
 $$
  \mbox{$
 V_+= \bigl( \frac{1+\a}{p-1}\bigr)^{\frac 1{p-1}}=
 5 \quad \mbox{for \, $p=2$ \, and \, $\a=4$}.
  $}
 $$

%%%%%%%%%%%%%%%%%%%%%%%%%%%%%%%%%%%%%%%%%%%%%%%%%%%%%%%%%%%%

\begin{figure}
%  \vskip -.3cm
\centering
\includegraphics[scale=0.65]{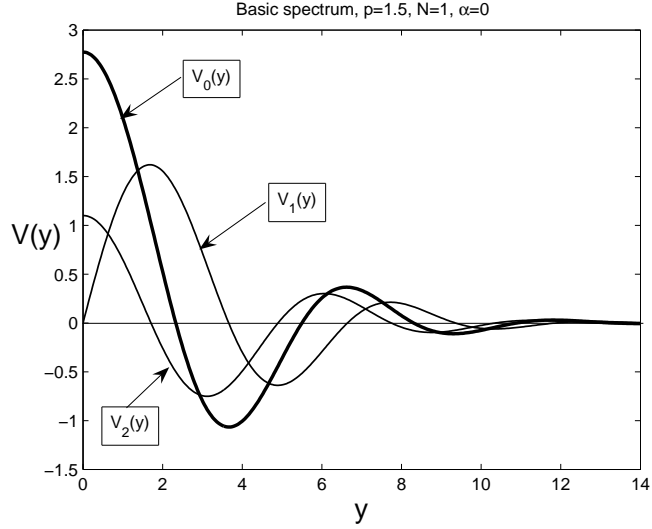}  %%%%%%%{4het5.eps}   %%%%%{4het.eps}  Old
\vskip -.5cm \caption{\small Three VSS profiles $V_0$, $V_1$, and
$V_2$ for $\a=0$; $p=1.5$ and $N=1$.}
   \vskip -.3cm
 \label{F30}
\end{figure}

For convenience, below we present the standard multiplicity result
for the autonomous case. In Figure \ref{F30}, we show first three
solutions $V_0(y)$, $V_1(y)$, and $V_2(y)$ of the problem
(\ref{4s})  with $p=1.5$  and
 %%\beq
 %%\label{aa1}
$ \a=0$.
 %% \eeq
  As was shown in \cite{GW2}, these $p$-bifurcation branches
  exhaust the whole set of VSS profiles in the case
  $\a=0$. %%%(\ref{aa1}).
 Notice that the mutual geometry of VSS profiles in the case
$\a=0$  %% (\ref{aa1})
   essentially differs
from that  in Figure \ref{F3}.
 %%% is impossible in the case
%%%%(\ref{aa1}).

%%%%%%%%%%%%%%%%%%%%%%%%%%%%%%%%%%%%%%%%%%%%%%%%%%%%%%%%%%%%

%%%%%%%%%%%%%%%%%%%%%%%%%%%%%%%%%%%%%%%%%%%%%%%%%%%%%%%%%

\begin{figure}
%  \vskip -.3cm
\centering
\includegraphics[scale=0.65]{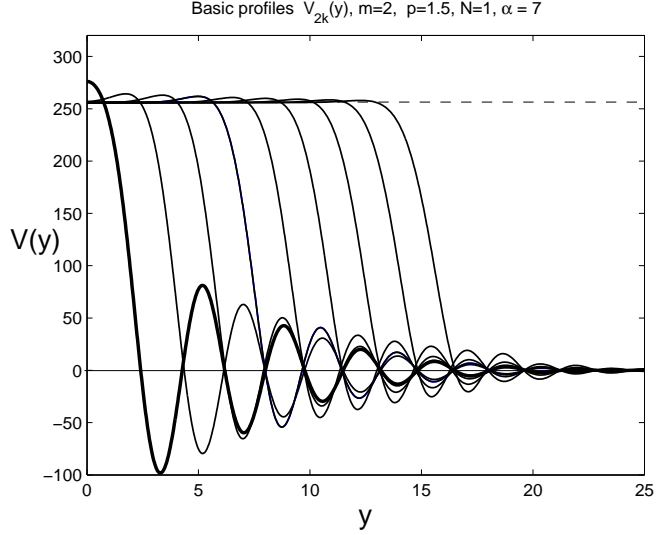}  %%%%%%%{4het5.eps}   %%%%%{4het.eps}  Old
\vskip -.5cm \caption{\small Several VSS profiles $V_{2k}(y)$ from
 satisfying (\ref{4s}) for $p=1.5$ and
$\a =7$.}
   \vskip -.3cm
 \label{F4}
\end{figure}

\begin{figure}
%  \vskip -.3cm
\centering
\includegraphics[scale=0.65]{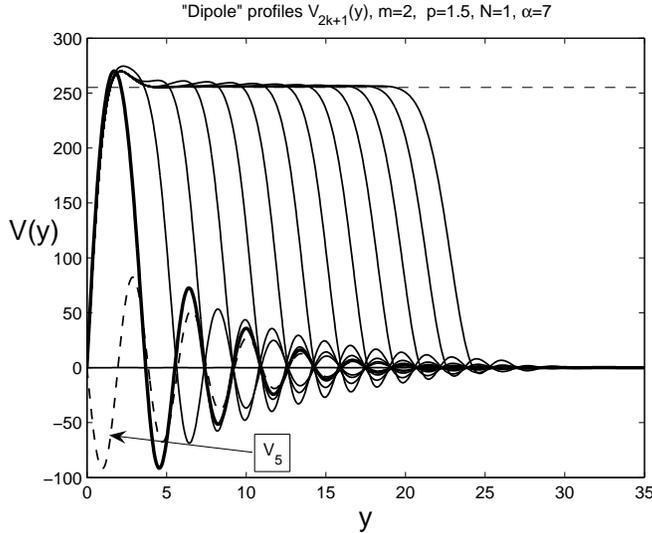}  %%%%%%%{4het5.eps}   %%%%%{4het.eps}  Old
\vskip -.5cm \caption{\small Several odd VSS profiles
$V_{2k+1}(y)$  satisfying (\ref{4s}) for $p=1.5$ and $\a =7$.}
   \vskip -.3cm
 \label{F5}
\end{figure}

 Let us more clearly show that, for $\a>0$, there occur other
  bifurcation-branching phenomena.
In Figures \ref{F4} and \ref{F5}, we again consider the case
$p=1.5$, $N=1$, and $\alpha=7$.
%%It turns out that, for such
For this sufficiently large $\a$, according to (\ref{al11}),
%%several saddle-node $\alpha$
 %%%bifurcations occur for smaller $\a < 7$, so
 we observe  several profiles from the basic  family
  $\{V_{2k}\}$ of even functions that are now concentrated about
  the constant equilibrium.
 %%splitting of basic lines.
  Namely, in Figure \ref{F4}, we show
nine  profiles $\{V_{2k}\}$.
%% that appeared from the basic
%%profiles $V_0^{(0)}$ by a cascade of similar bifurcations.
 Figure
\ref{F5} demonstrates
%%% even more complicated splitting and at least
thirteen odd profiles from the family
 $\{V_{2k+1}\}$.
%%appeared from the basic ``dipole-like" profile
 %%%$V_1^{(0)}$.
  In this case, instead of symmetry conditions
(\ref{Symm1}), we take the anti-symmetry ones (\ref{Asymm1}).
 Then the solutions $V(y)$ are continued in the odd manner,
 $V(-y) \equiv -V(-y)$ for $y<0$.

%%%%%%%%%%%%%%%%%%%%%%%%%%%%%%%%%%%%%%%%%%%%%%%%%%%%%%%%%%%%

%%%%%%%%%%%%%%%%%%%%%%%%%%%%%%%%%%%%%%%%%%%%%%%%%%%%%%%%%%%%

%%We do not have any analytical evidence on such saddle-node
%%%bifurcations.

%%%More numerics here ???????

\begin{figure}
%  \vskip -.3cm
\centering
\includegraphics[scale=0.65]{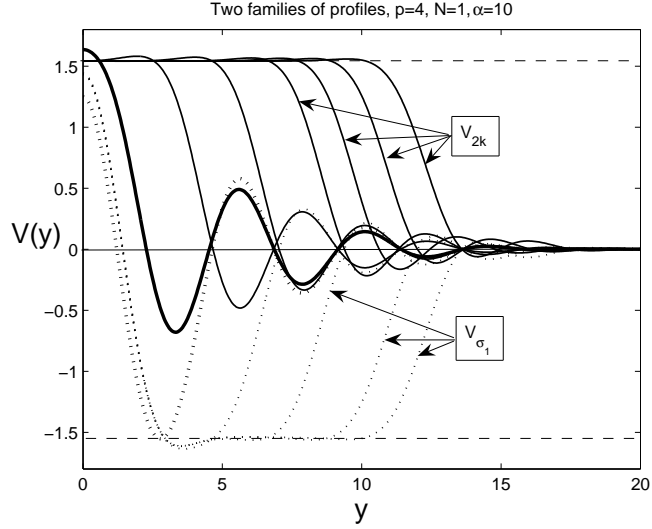}  %%%%%%%{4het5.eps}   %%%%%{4het.eps}  Old
\vskip -.5cm \caption{\small VSS profiles from $\{V_{2k}\}$ and
$\{V_{\s_1}\}$ for $p=4$, $N=1$, and $\a =10$.}
   \vskip -.3cm
 \label{F15}
\end{figure}

Before performing some easy preliminary estimates, we show in
Figure \ref{F15} the case $p=4$ and $\alpha=10$, where a similar
set of seven VSS patterns from the even family $\{V_{2k}\}$
 are presented.
 In addition, we see here
 another new family
  \beq
  \label{nn1}
  \{V_{\s_1}\}, \quad \mbox{with multiindex}
   \quad \s_1=\{+0,1,-k\},
    \eeq
where $0,1$, and $k$ in $\s_1$ stand for an ``effective number of
intersections" of the profiles with three consecutive equilibria
 $V_+,0$, and $-V_+$. We will discuss this new family later on.

%%from $\{V_2^{(k)}\}$  are presented. The profile with the question
%%mark ``?" shows that branching within such families can be more
%%complicated.
 %%% than as explained above.

 %%, together with two (dashed
%%lines) patterns from new families to be discussed later on.

%%%%%%%%%%%%%%%%%%%%%%%%%%%%%%%%%%%%%%%%%%%%%%%%%%%%%%%%%%%%

\begin{figure}
%  \vskip -.3cm
\centering
\includegraphics[scale=0.65]{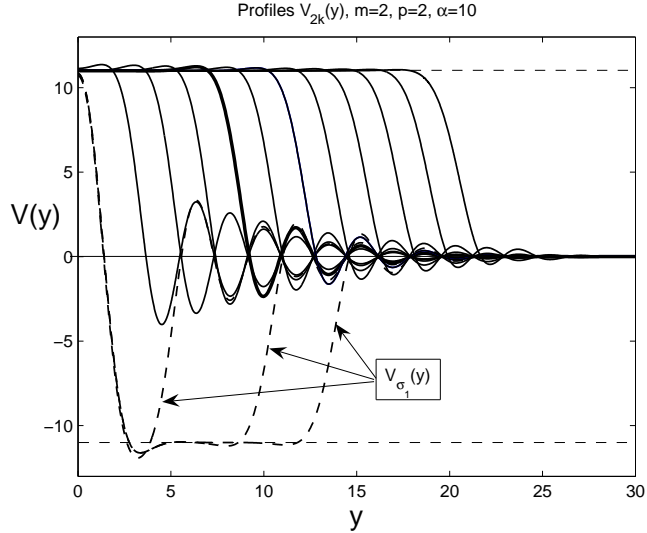}  %%%%%%%{4het5.eps}   %%%%%{4het.eps}  Old
\vskip -.5cm \caption{\small Profiles from the family $\{V_{2k}\}$
for $m=p=2$, $N=1$, and $\a =10$.}
   \vskip -.3cm
 \label{F66}
\end{figure}

In Figure \ref{F66}, we show the basic profiles $V_{2k}(y)$ and
others for $m=p=2$ and $\a=10$. This computation is important to
treat the results in Figure \ref{FF3} on the $\a$-branches.

%%%%%%%%%%%%%%%%%%%%%%%%%%%%%%%%%%%%%%%%%%%%%%%%%%%%%%%%%%%%

\begin{figure}
%\vskip -.3cm
\centering \subfigure[profiles $V_{2k}$]{
\includegraphics[scale=0.5]{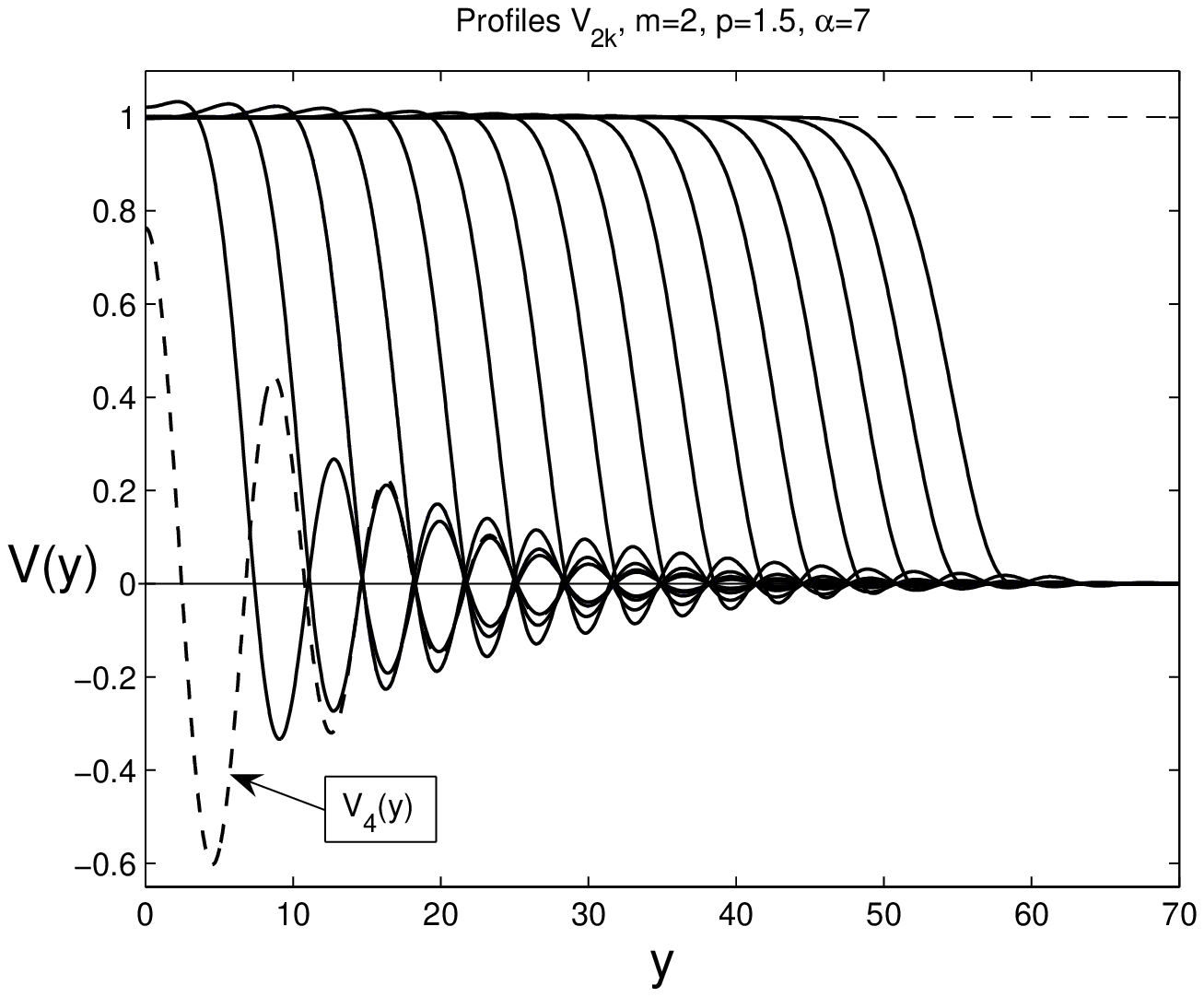} %Old: no N
} \subfigure[profiles $V_{\s_1}$]{
\includegraphics[scale=0.5]{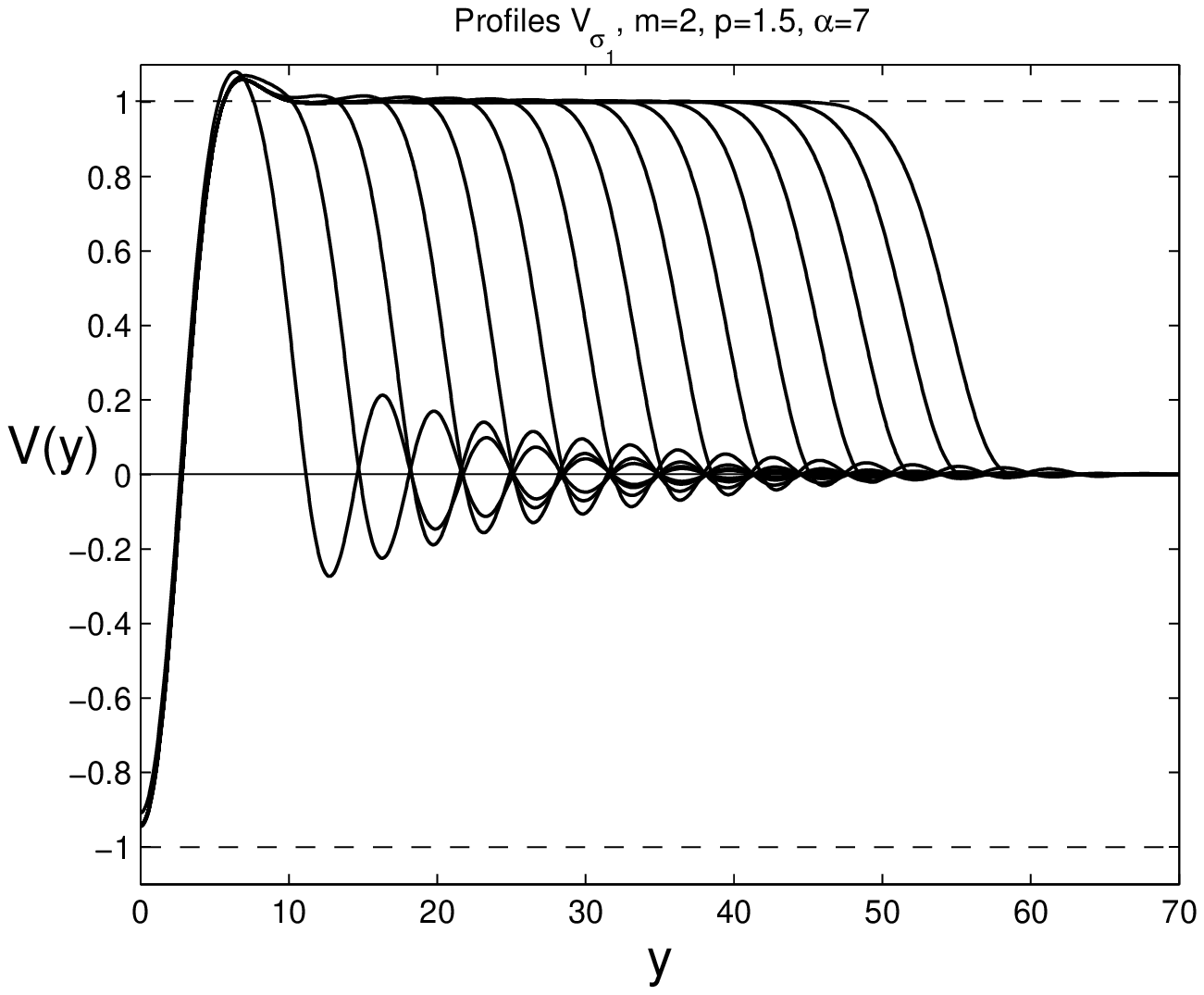}
}
 \vskip -.4cm
\caption{\rm\small Two different families  of VSS profiles for
$m=2$, $p=1.5$, $N=1$, and $\a=7$.}
  \vskip -.1cm
 %of the
 %ODE (\ref{pp.5}) for $m=-0.75$ (a) and $m=-0.9$ (b).}
% for
 %$n=1$ near non-oscillatory (a) and oscillatory (b) interfaces.}
 \label{FF33}
\end{figure}

\begin{figure}
%  \vskip -.3cm
\centering
\includegraphics[scale=0.65]{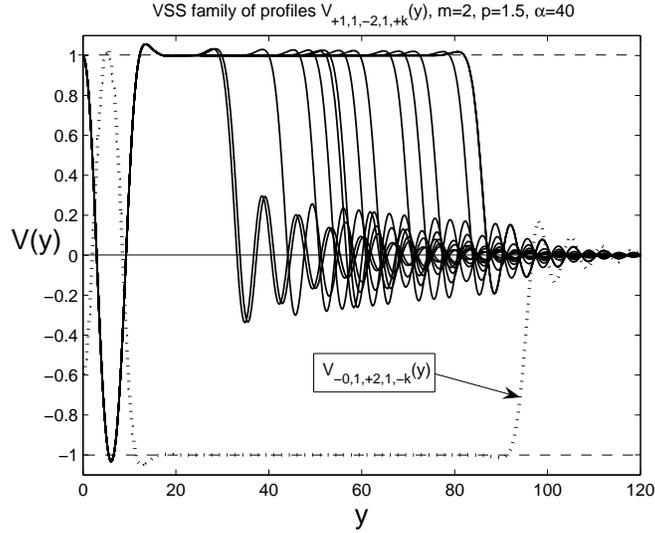}  %%%%%%%{4het5.eps}   %%%%%{4het.eps}  Old
\vskip -.5cm \caption{\small Some profiles from the family
$\{V_{\s_2}\}$ for $m=2$, $p=1.5$, and $\a =40$.}
   \vskip -.3cm
 \label{FF33N}
\end{figure}

For future convenience, some of the results are natural to present
for the rescaled equation by using the scaling
 \beq
 \label{ss2}
 V \mapsto C V, \quad y \mapsto c y, \quad \mbox{where}
  \quad C^{p-1}=\b, \,\,\, c= C^{\frac{p-1}{2m}}= \b^{\frac 1{2m}},
 \eeq
 that gives the following rescaled ODE:
  \beq
  \label{ss3}
   \mbox{$
 %%% {\bf D}_\infty (V) \equiv
  (-1)^{m+1} V^{(2m)} +  \frac 1 \b \,  \frac 1{2m} \, y
  V'  +  V - |V|^{p-1}V=0.
   $}
   \eeq
The constant equilibria are now fixed and are independent of the
parameters,
 \beq
 \label{eq11}
 V_\pm = \pm 1 \quad (\mbox{and} \,\,\, V=0),
 \eeq
 so we do not observe such huge VSS profiles as, for instance, in Figure
 \ref{F4}. The corresponding rescaled families $\{V_{2k}\}$ and
$\{V_{\s_1}\}$ for $m=2$, $p=1.5$, and $\a=7$ are shown in Figure
\ref{FF33}(a) and (b). Figure \ref{FF33N} shows another more
complicated family of even VSS profile
 $$
 V_{\s_2}(y), \quad \mbox{where}
  \quad \s_2=\{+1,1,-2,1,+k\}.
  $$
By boldface dotted line we denote therein  a different profiles
from the family $\{V_{\s_3}\}$, with a distinct multiindex $
\s_3=\{-0,1,+2,1,-k\}$. It seems that there are other families of
VSS profiles with more involved geometric structures that are very
difficult to catch numerically, so we stop the discussion at this
moment.

We expect that all the new profiles from families $\{V_{\s_1}\}$,
$\{V_{\s_2}\}$, etc.,
 are originated at  saddle-node bifurcations at some
 $\a_{\s_j}>0$ and cannot be extended to $\a=0$. Recall that such VSS
 solutions were not detected in the autonomous case $\a=0$,
 \cite{GW2}. Numerics with $\Delta \a=10^{-3}$ show some
 evidence concerning this though  are also rather unstable demonstrating
 typical jumps between many neighbouring $\a$-branches.
 We believe  that such  stable jumps definitely show
 existence of saddle-node points that are difficult to detect
 numerically.

%%%%%%%%%%%%%%%%%%%%%%%%%%%%%%%%%%%%%%%%%%%%%%%%%%
%%%%%%%%%%%%%%%%%%%%%%%%%%%%%%%%%%%%%%%%%%%%%%%%%%%%%%%%%%%%%%%%%%%%%%%%%

%%%%%%%%%%%%%%%%%%%%%%%%%%%%%%%%%%%%%%%%%%%%%%%%%%%%%%%%%%%%

%%\smallskip

%%\noi{\bf On some estimate: on countable number in each
%%$\a$-splitting.} First, let us estimate a total number of VSS
%%profiles from the typical families showed in Figures \ref{F4} and
%%\ref{F5}. According to this Figures, we assume that each $k$th
%%profiles $V^{(0,1)}_k(y)$ has the semi-width $k l$ where $l>0$ is
%%an approximate width of each layer added in the transition $k
%%\mapsto k+1$ for $k \gg 1$. Then multiplying the ODE (\ref{1.7R})
%%($m=2$, $N=1$) by $V$ in $L^2(\re)$ and integrating by parts
%%yields
%% $$
%% \mbox{$
%% \|V''\|_2^2+  \rho \|V\|^2_2 - \|V\|_{p+1}^{p+1}=0, \quad
%% \mbox{where}
%% \quad \rho= \frac{1+\a}{p-1}- \frac N{2m}.
%%  $}
%%  $$
%%  Since the first term is negligible in view of the almost ``flat"
%%shape of $V_k$ for $k \gg 1$ and estimating the other two by using
%%the given layer structure of patterns, we obtain the approximate
%%inequality
 %%$$
 %%\rho V_*^2 k l \approx V_*^{p+1} k l,
 %%$$
%%so this does not characterize a possible maximal number $k$, i.e.,
%%the cardinal number of such VSS sets. Actually, the general
%%structure of the ODE (\ref{1.7R}) indicates that such sets can be
%%infinite (countable).

%%We were not able to detect other families of solutions of
%%(\ref{ss3}) rather than the above basic family $\{V_k\}$ and
%%$\{V_{\s_1}\}$  as explained in Figure \ref{FF33}(b).

\smallskip

%%%%%%%%%%%%%%%%%%%%%%%%%%%%%%%%%%%%%%%%%%%%%%%
\noi {\bf On countable sets via shooting arguments.} As an
illustration, consider the ODE in (\ref{4s}) (i.e., for $m=2$ and
$N=1$) in the symmetric case (\ref{Symm1}), so we study $V_{2k}$
and other even profiles. The shooting argument justifies existence
of a connection of the 2D bundle near the origin, which is
characterized by two symmetry conditions (\ref{Symm1}), so that
the bundle has two parameters
 \beq
 \label{ff1}
 f(0) \quad \mbox{and} \quad f''(0).
 \eeq
 We then need to intersect this bundle
 with the exponential bundle at $y= + \infty$, where the linearized ODE
 admits the following WKBJ-type asymptotics (we omit power-like
 multipliers that are not important at this stage):
  $$
   \mbox{$
  V^{(4)}= - \frac 1{2m} \, V' y +...
 \quad \Longrightarrow \quad V(y) \sim {\mathrm e}^{a y^{4/3}}, \quad
 \mbox{with} \,\,\,
 a^3 = \frac 14\big(\frac 34\big)^3.
  $}
  $$
  Therefore, there exists a 1D unstable bundle
with positive $a_0= \frac 34 \, 4^{-1/3}$,
   and a 2D stable oscillatory bundle with
   $$
    \mbox{$
   a_\pm= a_0\big(-\frac 12 \pm {\mathrm i}\frac{\sqrt 3}2\big) \quad\,\,
   \big({\rm Re} \, a_\pm < 0 \big).
    $}
    $$
Thus, we observe here a typical difficult two-parameter shooting
problem of 2D$\to$2D connection of two   manifolds. In this case,
in the above shooting approach, we arrive at two equations for two
parameters given in (\ref{ff1}). In the analytic case $p=3,5,...$,
where the dependence  on parameters is  assumed to be analytic, we
conclude that {\em the set of VSS profiles is at most countable.}
 Moreover, for fixed $p<p_0$ and any $\a>0$, we expect a finite
(but, possibly, arbitrarily large for some values) number of VSS
profiles according to the above $p$- and $\a$-bifurcations and
extra possible $\a$-saddle-node bifurcations that much less are
known about.

\smallskip

\begin{figure}
%  \vskip -.3cm
\centering
\includegraphics[scale=0.75]{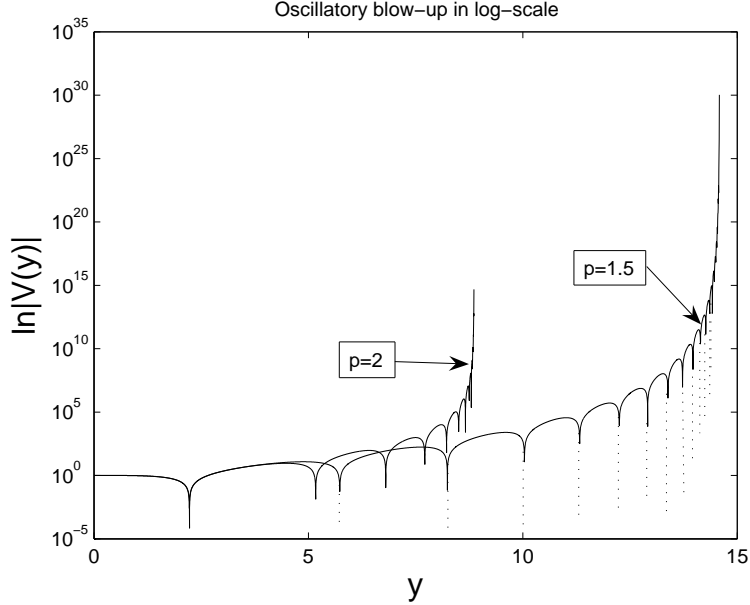}  %%%%%%%{4het5.eps}   %%%%%{4het.eps}  Old
\vskip -.5cm \caption{\small
 Oscillatory blow-up for the ODE (\ref{Bl1})
 for $p= \frac 32$ and $p=2$.}
   \vskip -.3cm
 \label{FBl1}
\end{figure}

%%%%%%%%%%%%%%%%%%%%%%%%%%%%%%%%%%%%%%%%%%%%%%%%%%%%%%%%%%%%%%%%%%%%
\noi{\bf On oscillatory blow-up in the ODE.} Here we briefly
discuss an important aspect of ODEs such as (\ref{1.7R}), $N=1$
(or in radial geometry). Namely, these admit a strong (nonlinear)
unstability of orbits $V(y)$  that  is associated not with the
exponential growth in the linear setting, but with the blow-up due
to the main two terms in the ODE in (\ref{4s}). For $m=2$ these
are
  \beq
  \label{Bl1}
 V^{(4)} = -|V|^{p-1} V,
   \eeq
 where we omit all linear terms that are negligible as $|V| \to +
 \infty$. Obviously, (\ref{Bl1}) does not admit blow-up solutions
 of constant sign, so we need to describe the {\em oscillatory
 blow-up}. In Figure \ref{FBl1}, we present typical blow-up behaviour
 of solutions $\ln|V(y)|$
 of (\ref{Bl1}) for $y>0$ with Cauchy data
 $$
 V(0)=1, \quad V'(0)=V''(0)=V'''(0)=0.
 $$
Solutions blow-up at some finite $y_0>0$ with definite oscillatory
behaviour as $y \to y_0$.

%%%%%%%%%%%%%%%%%%%%%%%%%%%%%%%%%%%%%%%%%%%%%%%%%%%%%%%%%%%%
%%%%%%%%%%%%%%%%%%%%%%%%%%%%%%%%%%%%%%%%%%%%%%%%%%%%%%%%%%%%%%%%%%%%%%%%

In order to detect this oscillatory behaviour,
 by shifting $y_0-y \mapsto -y$, we have that blow-up
 happens at $y=0$, i.e., as $y \to 0^-$, and introduce the
 {\em oscillatory component} $\varphi(s)$ as follows:
  \beq
  \label{Bl2}
   \mbox{$
 V(y) = (-y)^\mu \varphi(s), \quad \mbox{where} \quad s= \ln (-y)
  \quad \mbox{and} \quad \mu= - \frac 4{p-1}<0.
   $}
  \eeq
Then $\varphi(s)$ solves the following autonomous ODE in $\re$:
 \beq
 \label{Bl3}
  \begin{matrix}
 \varphi^{(4)} + 2(2\mu-3) \varphi''' + (6 \mu^2- 18 \mu +11) \varphi''
 + 2(2 \mu^3
   - 9 \mu^2  \smallskip\smallskip\cr + \, 11 \mu -3) \varphi'
 + \mu(\mu-1)(\mu-2)(\mu-3) \varphi = -|\varphi|^{p-1} \varphi.
  \end{matrix}
   \eeq
 Obviously, this ODE also admits blow-up singularities as
 (\ref{Bl1}), but we are interested in bounded or slower growing  orbits that can be extended up to
 the blow-up point $s=-\infty$.

 In addition, using a standard min-max scenario of orbits behaviour on parameters, we
 expect a bounded (hopefully, periodic) or a spiral-type solution $\varphi_*(s)$ in $\re$ to exist,
 which, according to (\ref{Bl2}), will describe a generic
 structure of blow-up singularities. The blow-up instabilities make it very
 difficult to see necessary orbits of equation
 (\ref{Bl3}) even numerically. A logarithmic-type distribution of
 zeros close to blow-up points as in (\ref{Bl2}) is seen in Figure
 \ref{FBl1}.

  What is crucial in this
 case, is that the whole unstable bundle is {\em two-dimensional}
 incorporating the parameter $y_0$ as the position of the blow-up
 point (i.e., $y \mapsto y_0-y$ in (\ref{Bl2})) and the parameter
 $s_0$ of translation, $\varphi_*(s)\mapsto \varphi_*(s+s_0)$.
Therefore, our VSS profiles are those that  have nothing to do
with this singular blow-up bundle, so we again observe two
equations for two parameters in (\ref{ff1}).

%%%%%%%%%%%%%%%%%%%%%%%%%%%%%%%%%%%%%%%%%%%%%%%%%%%
\section{Passing to the limit $\a \to \infty$}
 \label{SectInf}

Without loss of generality, we consider the 1D case in
(\ref{1.7R}) and bear in mind the first symmetric similarity
profile
 $V_0(y)$ for the
ODE
 \beq
 \label{ss1}
  \mbox{$
 (-1)^{m+1} V^{(2m)}  + \frac 1{2m} \, y V' + \b V - |V|^{p-1}V=0, \quad
 \mbox{where} \quad
 \b= \frac{1+\a}{p-1} \to + \infty
  $}
  \eeq
  as $\a \to + \infty$.
  For $\b \gg 1$,
the first natural  step is to use the scaling (\ref{ss2})
 %%\beq
%% \label{ss2}
%% V \mapsto C V, \quad y \mapsto c y, \quad \mbox{where}
%%  \quad C^{p-1}=\b, \,\,\, c= C^{\frac{p-1}{2m}},
%% \eeq
 and to write the resulting equation in the   perturbed
 form,
  \beq
  \label{ss3N}
   \mbox{$
  {\bf D}_\infty (V) \equiv
  (-1)^{m+1} V^{(2m)}  +  V - |V|^{p-1}V= -  \frac 1 \b \,  \frac 1{2m} \, y
  V'.
   $}
   \eeq
 Hence,  the right-hand side becomes negligible as $\b \to \infty$ on
 smooth solutions.

Formally passing to the limit $\b \to \infty$ in (\ref{ss3N}), we
obtain the {\em limit equation} \beq
  \label{ss3Nl}
   \mbox{$
  {\bf D}_\infty (W) \equiv
  (-1)^{m+1} W^{(2m)}  +  W - |W|^{p-1}W=0 \quad \mbox{for}
  \quad y>0, \quad W(+\infty)=0.
   $}
   \eeq
Let us show that a nontrivial solution of (\ref{ss3Nl}) does not
exist. For instance, for $m=2$,
 \beq
 \label{e1}
  \begin{matrix}
 -W^{(4)}+W-|W|^{p-1}W=0, \quad \mbox{so that, for $y \gg 1$,}
 \smallskip\smallskip\\
 %%\eeq
 %%% so that, for $y \gg 1$,
 %% we have
  %%%$$
  W^{(4)}=W+... \quad \Longrightarrow \quad \exists \,\, 1D \,\,
  \mbox{bundle}, \,\,\, W(y) \sim C {\mathrm e}^{-y}+... \, .
 \end{matrix}
 \eeq
  %%$$
  This 1D is not enough to satisfy two symmetry conditions
  (\ref{Asymm1}).
  Similarly, for $m=3$,
\beq
 \label{e2}
  \begin{matrix}
 W^{(6)}+W-|W|^{p-1}W=0, \quad \mbox{and for $y \gg 1$,}
 \smallskip\smallskip\qquad\quad\\
 %%\eeq
 %%%% and for $y \gg 1$,
 %% $$
  W^{(6)}=-W+... \quad \Longrightarrow \quad \exists \,\, 2D \,\,
  \mbox{bundle with char. values $\l_\pm=-\frac \pi 6 \pm {\rm i} \, \frac {\sqrt 3}2$}.
  \qquad\quad
  %% \,\,\, V(y) \sim C {\mathrm e}^{-y}+... \, .
   \end{matrix}
   \eeq
  The 2D bundle is not sufficient to satisfy three symmetry
  conditions at the origin,
 $$
 W'(0)=W'''(0)=W^{(5)}(0)=0.
  $$
  The same lack of dimension happens for any $m \ge 1$.
In other words, nonexistence for
 (\ref{ss3Nl}) is a typical property.

\begin{figure}
%  \vskip -.3cm
\centering
\includegraphics[scale=0.65]{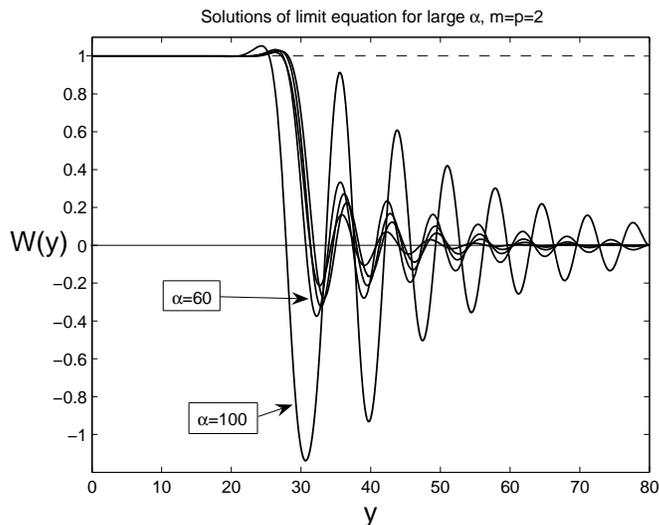}
\vskip -.5cm \caption{\small
 Towards nonexistence for the ODE (\ref{e1}), $m=p=2$: solutions of
(\ref{ss3N}) get very oscillatory for $\a \gg 1$.}
   \vskip -.3cm
 \label{FN1}
\end{figure}

Such nonexistence is easily checked numerically. In Figure
\ref{FN1}, we show for $m=p=2$ that as $\a$ increases the VSS
profiles with a fixed length of the flat part (such solutions
exist) become more and more oscillatory and ceases  to exist for
$\a=+\infty$.

%%%%%%%%%%%%%%%%%%%%%%%%%%%%%%%%%%%%%%%%%%%%%%%%%%%%%%%%%%%%
%%%%%%%%%%%%%%%%%%%%%%%%%%%%%%%%%%%%%%%%%%%%%%%%%%%%%%%%%%%%%%%%%%%%%%%%

\begin{figure}
%  \vskip -.3cm
\centering
\includegraphics[scale=0.65]{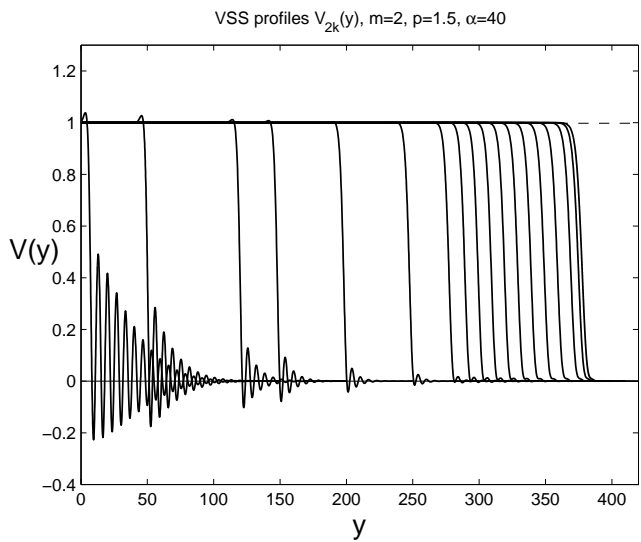}
\vskip -.5cm \caption{\small
 Some VSS profiles satisfying
(\ref{ss3N}) for $m=2$, $p=1.5$, and $\a = 40$.}
   \vskip -.3cm
 \label{FN2}
\end{figure}

On the other hand, profiles with huge  flat non-oscillatory parts
get very wide for $\a \gg 1$. For instance, Figure \ref{FN2} shows
for $p=1.5$ that for $\a=40$, the last less oscillatory VSS
profile $V_{2k}(y)$ has the flat part of the length $\sim 380$,
while the first profiles are strongly oscillatory (and do not
admit passage to the limit $\a \to \infty$). Notice that according
to bifurcation exponents (\ref{al11}), for $\a=40$, we expect, at
least,
 $$
 \mbox{164 profiles $V_{2k}$ and 164 profiles $V_{\s_1}$,}
 $$
 i.e., 328 VSS profiles overall (and more in view of
possible saddle-node bifurcations beforehand, at some $\a_{\rm
s-n} \in (0, 40)$).

%%%%%%%%%%%%%%%%%%%%%%%%%%%%%%%%%%%%%%%%%%%%%%%%%%%%%%%%%%%%
%%%%%%%%%%%%%%%%%%%%%%%%%%%%%%%%%%%%%%%%%%%%%%%%%%%%%%%%%%%%%%%%%%%%%%%%

%%%Nonexistence for $\a=+\infty$ is also easily seen analytically by
%%%Pohozaev's Identity \cite{Poh00}.

Therefore, as $\a \to +\infty$, one can that the VSS profiles get
wider indefinitely in the sense that, in the rescaled variables,
uniformly on bounded intervals,
 \beq
 \label{ff1s}
 V_0(y) \to 1 \quad \mbox{as} \quad \a \to + \infty,
  \eeq
  and the divergence is sufficiently (exponentially) fast.
 %% The rate of this divergence of the semi-width of the profile is
 %% unknown and can be obtained by matching for the ODE
 %% (\ref{ss3N}).

As a simple illustration of this behaviour presented in Figure
\ref{FN2}, we give a nonexistence result concerning a step-like
profile for the limit equation (\ref{e1}). Assume that, unlike
(\ref{ff1s}), there exists a nontrivial finite limit as $\a \to
\infty$, so (\ref{e1}) admits solution $W(y)$ with exponential
decay such that
 \beq
 \label{WW1}
 W(y) \approx 1 \quad \mbox{for}
 \quad y \in[0,L], \,\,\, \mbox{with some} \,\,\, L \gg 1,
 \eeq
 and that the internal layer at $y \sim L$ does not depend essentially on $L$
 (recall that (\ref{e1}) is an autonomous ODE), and, in a natural sense,
  \beq
  \label{WW2}
   \mbox{$
   %%\int (W'')^2, \,\,
    \int_L^\infty W^2 \ll L \quad \mbox{and}
    \quad \int_L^\infty |W|^{p+1} \ll L .
   $}
  \eeq
The following nonexistence result has a ``conditional" nature,
where we prohibit certain solutions of  a specific spatial shape:

\begin{proposition}
 \label{Pr.GG1}
 The ODE $(\ref{e1})$
 does not admit a nontrivial solution with exponential decay at
 infinity
  satisfying $(\ref{WW1})$ and
 $(\ref{WW2})$.

  \end{proposition}

  \noi{\em Proof.} We use two identities obtained by multiplying
  (\ref{e1}) in $L^2$ by $W$ and $y W'$ (Pohozaev's multiplier in
  elliptic theory, \cite{Poh00}) to get
   \beq
   \label{PP1}
   \left\{
    \begin{matrix}
    - \int (W'')^2 + \int W^2 - \int|W|^{p+1}=0,\qquad\quad
    \smallskip\smallskip \\
 -\frac 32  \int (W'')^2 - \frac 12   \int W^2 +\frac 1{p+1}
 \int|W|^{p+1}=0.
  \end{matrix}
  \right.
  \eeq
  Substituting $\int(W'')^2$ from the first into the second
  identity yields
 \beq
 \label{PP2}
 \mbox{$
 \int W^2= \g_0
 \int|W|^{p+1},
  \quad \mbox{where} \quad
 \g_0= \frac {3p+5}{4(p+1)}.
  $}
  \eeq
  %%%from whence follows $W=0$. $\qed$
 Since $\g_0<1$ for any $p>1$, this contradicts hypotheses (\ref{WW1}),
(\ref{WW2}) under which both integrals in (\ref{PP2}) are close to
$L$. $\qed$

%%%%%%%%%%%%%%%%%%%%%%%%%%%%%%%%%%%%%%%%%%%%%%%%%%%%%%%%%%%%%%%%%%%%%%
%%%%%%%%%%%%%%%%%%%%%%%%%%%%%%%%%%%%%%%%%%%%%%%%%%%%%%%%%%%%%%%%%%%%%%%%%%%%%5
\section{On some  VSSs results for the non-monotone
absorption}
%%%$-|u|^p$}
 %%% $p \ge p_0$}
 \label{SectNon}

In this short section, we orient  our nonexistence business to the
following PDE:
 \beq
 \label{31}
 u_t = - (-\D)^m u - t^\a |u|^p,
  \eeq
where, unlike (\ref{1.1R}), we take the non-monotone absorption
nonlinearity $|u|^p$.

%%%%%%%%%%%%%%%%%%%%%%%%%%%%%%%%%%%%%%%%%%%%%%
\subsection{On local existence close to bifurcation points}

 Local existence of similarity solutions for
$p < p_0$ for equation (\ref{1.7R}),  where
 \beq
  \label{rep1}
 |V|^{p-1}V \mapsto |V|^p,
 \eeq
 can be established as above by bifurcation theory. Then,
 for several types of ``positively dominant" VSS profiles $V(y)$ as
in Figures \ref{F3}, \ref{F4}, \ref{F5}, \ref{F66}, and
\ref{FF33}(a), their structure changes slightly after the
replacement (\ref{rep1}).

\begin{figure}
%  \vskip -.3cm
\centering
\includegraphics[scale=0.65]{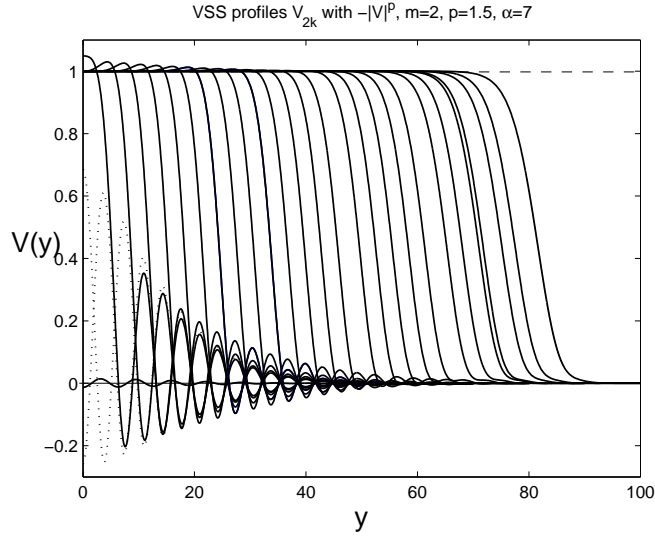}
\vskip -.5cm \caption{\small Twenty nine VSS solutions $V_{2k}(y)$
of
 (\ref{N4s}) for $m=2$, $p=1.5$, $\a=7$.}
   \vskip -.3cm
 \label{Fmod1}
\end{figure}

\begin{figure}
%  \vskip -.3cm
\centering
\includegraphics[scale=0.65]{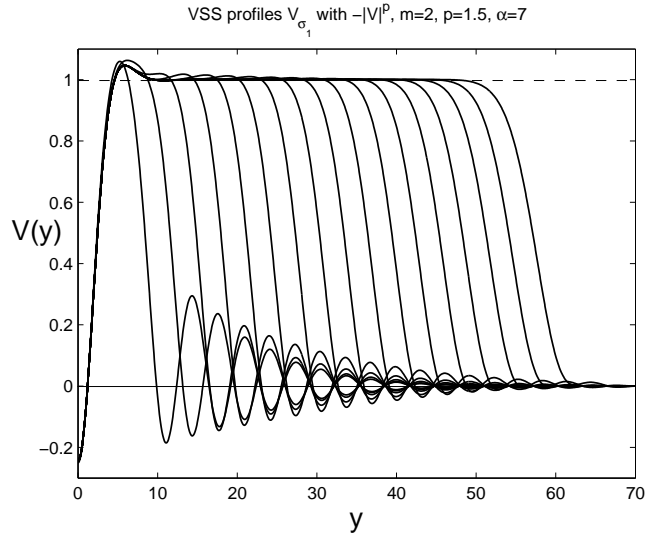}
\vskip -.5cm \caption{\small Seventeen VSS solutions $V_{\s_1}(y)$
of
 (\ref{N4s}) for $m=2$, $p=1.5$, and $\a=7$.}
   \vskip -.3cm
 \label{Fmod2}
\end{figure}

For instance, in Figure \ref{Fmod1}, we present even VSS profiles
for the ODE
 \beq
 \label{N4s}
  \mbox{$
 -V^{(4)} + \frac 14 \, y V' + \b V -|V|^{p}=0 \quad \mbox{in}
 \quad \re \quad \big(\b= \frac{1+\a}{p-1}\big),
 $}
  \eeq
in the case $p=1.5$, $\a=7$, where most of the solutions  are very
close to those  in Figure \ref{FF33}(a). Excluding first
sufficiently small profiles, all the others are clearly positively
dominant (the last ones look ``almost" positive; careful checking
shows that all of them remain oscillatory). Small solutions
(dotted lines as well as that of order $\sim 2\cdot 10^{-2}$ given
by the solid line) become essentially different. Instead of
(\ref{eq11}),  equation (\ref{N4s}) admits the single constant
non-zero equilibrium
 \beq
 \label{nn11}
 V_+=+1 \quad (\mbox{and} \,\,\, V_0=0),
  \eeq
 so that VSS profiles $V_{\s_1}$ such as in Figure \ref{F15} and \ref{F66} are impossible.
 Since (\ref{nn11}) means that for the nonlinearity $-|V|^p$, two
 equilibria $-1$ and 0 just coincide (in fact, identified), we
 obtain another but similar family of profiles again denoted by
 $\{V_{\s_1}\}$ that are shown in Figure \ref{Fmod2}. For large
 $\a=40$, step-like solutions of (\ref{N4s}) are practically
 indistinguishable from
  %%%%coincide with
 those in Figure \ref{FN2}, so that this behaviour as $\a \to +
 \infty$ suits both monotone and non-monotone absorption terms.

%%%%%%%%%%%%%%%%%%%%%%%%%%%%%%%%%%%%%%%%%%%%%%%%%%%%%%%%%%%%
%%%%%%%%%%%%%%%%%%%%%%%%%%%%%%%%%%%%%%%%%%%%%%%%%%%%%%%%%%%%%%%%%%%%%%%%

%%%%%%%%%%%%%%%%%%%%%%%%%%%%%%%%%%%%%%%%%%%%%%%%%%%%%%%%%%%%
%%%%%%%%%%%%%%%%%%%%%%%%%%%%%%%%%%%%%%%%%%%%%%%%%%%%%%%%%%%%%%%%%%%%%%%%

On the contrary, nonexistence of VSSs for
 \beq
 \label{32}
  \mbox{$
 p \ge p_0= 1 + \frac{2m(1+\a)}N
 $}
 \eeq
can be proved much easier than for the original PDE (\ref{1.1R}).
This is explained in detail in \cite{GSVSS1} for $\a=0$, so we
present here a few nonexistence comments.

%%%%%%%%%%%%%%%%%%%%%%%%%%%%%%%%%%%%%%%%%%%%%%%%%%%%%%%%%%%%%%
\subsection{On nonexistence of similarity profiles}

The similarity solutions remain the same, (\ref{1.6R}), with the
elliptic equation (\ref{1.7R}), where we use the replacement
(\ref{rep1}).
 Therefore, integrating this equation over $\ren$ with the
 condition of exponential decay,
 %%% using
 %% (\ref{ExpC}),
  we obtain in particular that
  \beq
  \label{33}
  \mbox{$
  \int |V|^p=0 \quad \mbox{for}
  \quad p=p_0 \quad \Longrightarrow \quad V=0.
  $}
   \eeq
 Therefore, all the continuous $p$-branches do not cross
the vertical line $p=p_0$, so that in the parameter region $\{p
\ge  p_0\}$ there are no VSS profiles at all.

Obviously, for the original equation (\ref{1.7R}) this analysis
leads to
 \beq
  \label{33S}
  \mbox{$
  \int |V|^{p-1}V=0 \quad \mbox{for}
  \quad p=p_0,
  %%% \quad \Longrightarrow \quad V=0.
  $}
   \eeq
   which is  not that controversial.
Of course, this prohibits clearly positively dominant VSS profiles
such as those in Figure \ref{FF33}(a), but not obvious for
profiles having essential negative parts as in (b) therein.
Nevertheless, one can see that, for all these profiles,
(\ref{33S}) is valid, since these are also positively dominant.

%%%%%%%%%%%%%%%%%%%%%%%%%%%%%%%%%%%%%%%%%%%%
\subsection{Nonexistence for the PDE}

Using a slight modification of the statement and the proof in
\cite[\S~2]{GSVSS1} based on Pohozaev's nonlinear capacity
approach (see details and related references therein), we
formulate the following nonexistence result for (\ref{31}):

%%%%%%%%%%%%%%%%%%%%%%%%%%%%%%%%%%
\begin{theorem}
\label{Th.2}
 Let $\a>-1$ and
  \beq
  \label{ll12}
 \mbox{$
p \ge p_0= 1 + \frac{2m(1+\a)} N. $}
 \eeq
If a function $u(x,t) \in L^p_{\rm loc}(\ren \times \re_+
\setminus\{(0,0)\})$ satisfies $(\ref{31})$ and $(\ref{g15})$ in
the weak sense,
%% and
%%  \beq
%%   \label{g15}
%%u(x,0)=0 \quad \mbox{for all} \,\,\, x \not = 0,
%% \eeq
%%% in the sense of distributions,
 %%   in the sense of
%%%$(\ref{we1})$,
then $u=0$.
 \end{theorem}

In particular, this means that, for this nonlinearity, any
similarity problem such as (\ref{1.7R})
%%, (\ref{ExpC})
 does not
admit a nontrivial solution. Of course, this implies that no
other, non-similarity, VSSs exist. For the monotone nonlinearity
as in (\ref{1.7R}), the proof of nonexistence gets  more involved,
and can be done, with some technical changes,  along the lines in
\cite[\S~3.2]{GSVSS1}.
 %%%%  and is not presented here.

%%In this energy class, we are looking for solutions satisfying the
%%following  singular initial condition:

%%%%%%%%%%%%%%%%%%%%%%%%%%%%%%%%%%%
%%\subsection{Monotone nonlinearity}

%%We return to the original model (\ref{1.1R}), (\ref{al.1}). Here
%%we adapt the methods from \cite{GSVSS1}.

%%\com{AES: your turn, difficult ???}

\ms

{\bf Acknowledgement.}
 The author thanks A.E.~Shishkov for drawing his
 attention to time dependent absorption phenomena.
 %% and for
 %%%%useful discussions.
 %%%%%for  a stimulating discussion.
%%%%%%%%%%%%%%%%%%%%%%%%%%%%%%%%%%%%%%%%%%%%%%%%%%%%%%

%%%%%%%%%%%%%%%%%%%%%%%%%%%%%%%%%%%%%%%%%%%%%%%%%%%%%%%%%%%%%%%%%%%%%%%%%%%

%%%%%%%%%%%%%%%%%%%%%%%%%%%%%%%%%%%%%%%%%%%%%%%

\begin{thebibliography}{10}




%\bibitem %[JFT]
%{JFT} C. J. Amick and J. F. Toland {\em Homoclinic orbits in the
%dynamic phase space analogy of an elastic strut}, Euro. J. Appl.
%Math., {\bf 3} (1992), 97-114.


 \bibitem
{Rynn2}
 R.~Bari and B.~Rynne, {\em Solution curves and exact multiplicity results
 for $2m$th order boundary value problems},  J.~Math. Anal. Appl., {\bf 292} (2004), 17--22.


\bibitem %[Berger]
 {Berger}
  M.~Berger, {Nonlinearity and Functional Analysis}, Acad.
  Press, New York, 1977.





\bibitem %[BS]
{BS}
 M.S.~Birman and M.Z.~Solomjak, {Spectral Theory of Self-Adjoint Operators
 in Hilbert Space}, D. Reidel, Dordrecht/Tokyo, 1987.



 \bibitem{BrFr}
  H. Brezis and A. Friedman, {\em Nonlinear parabolic
  equations involving measures as initial conditions}, J. Math. pures
  appl., {\bf 62} (1983), 73--97.




\bibitem %[BPT]
{BPT} H.~Brezis, L.A.~Peletier, and D.~Terman, {\em A very
singular solution of the heat
  equation with absorption}, Arch. Rat. Mech. Anal., {\bf 95}
  (1986), 185--209.


\bibitem{BGW1}
C.J.~Budd, V.A.~Galaktionov, and J.F.~Williams, \emph{Self-similar
blow-up
  in higher-order semilinear parabolic equations}, SIAM J.~Appl. Math.,
 {\bf 64} (2004), 1775--1809.





\bibitem
{CodL} E.A. Coddington and N. Levinson, {\rm Theory of Ordinary
Differential Equations}, McGraw-Hill Book Company, Inc., New
York/London, 1955.



%%\bibitem
%%{KDal} Ju.L. Daleckii and  M.G. Krein, {\rm Stability of Solutions
%%of Differential Equations
%% in Banach Space},
%%Transl.  Math. Monographs, Vol.  {\bf 43}, Amer. Math. Soc.,
%%Providence, R.I., 1974.

\bibitem
{Deim}
K. Deimling, {\rm Nonlinear Functional Analysis}, Springer-Verlag,
Berlin/Tokyo, 1985.

 %%\bibitem %AUTO
%%{AUTO}
%%E.J. Doedel, A.R. Champneys, T.F.  Fairgrieve, Y.A. Kuznetsov,
%%  B. Sandstede, and  X.-J. Wang.  {\em AUTO97: Continuation and bifurcation software for %%ordinary
%%  differential equations}, Technical report, Department of Computer
%%  Science, Concordia University, Montreal, Canada.
%% Available at {\tt ftp://ftp.cs.concordia.ca/pub/doedel/auto}.
%\bibitem %[AUTO]
%{AUTO} E.J. Doedel, {\em AUTO, a program for the automatic
%continuation of autonomous systems}, Cong. Numer., {\bf 30}
%(1981), 265-384.



%\bibitem %[EGKP1]
%{EGKP1} Yu.V. Egorov,  V.A. Galaktionov, V.A. Kondratiev, and S.I.
%Pohozaev, {\em  On the necessary conditions of existence to a
%quasilinear inequality in the half-space,} Comptes Rendus Acad.
%Sci. Paris,  S\'erie I, {\bf 330} (2000), 93-98.



  \bibitem{Eg4}
Yu.V.~Egorov, V.A.~Galaktionov, V.A.~Kondratiev, and
S.I.~Pohozaev,
  {\em Asymptotic behaviour of global solutions to higher-order semilinear
  parabolic equations in the supercritical range}, Adv. Differ. Equat.,
{\bf 9} (2004), 1009--1038.





 %%\bibitem %[Eg4]
 %%{Eg4}
%% Yu.V. Egorov,  V.A. Galaktionov, V.A. Kondratiev, and S.I. Pohozaev, {\em
%% Asymptotic behaviour of global solutions to higher-order semilinear parabolic equations in
%% the supercritical range},  {Comptes Rendus
%%Acad. Sci. Paris, S\'erie I}, {\bf 335} (2002), 805-810 (full text
%%in {\tt http://www.maths.bath.ac.uk/ MATHEMATICS/preprints.html}).
% {\tt http://mip.ups-tlse.fr}).



\bibitem %[EidSys]
{EidSys} S.D. Eidelman, {Parabolic Systems,}  North-Holland Publ.
Comp., Amsterdam/London, 1969.


\bibitem
{Ellias} U.~Ellias, {\em Eigenvalue problems for the equation $L y
+ \l p(x) y=0$}, J. Differ. Equat., {\bf 29} (1978), 28--57.



\bibitem
 {EscKav} M.~Escobedo and O.~Kavian, {\em Variational problems
 related to self-similar solutions of the heat equation},
 Nonlinear Anal., TMA, {\bf 11} (1987), 1103--1133.

 %% \bibitem
%%  {EGW1}
%%J.D.~Evans, V.A.~Galaktionov, and J.F.~Williams, \emph{Blow-up and
%%global asymptotics of  the limit unstable Cahn-Hilliard equation},
%%SIAM J. Math. Anal., {\bf 38} (2006), 64--102.




%%\bibitem
%% {Fedor} M.V. Fedoryuk, {\em  Singularities of the kernels of Fourier integral
%%operators and the asymptotic behaviour of the solution of the
%%mixed problem}, Russian Math. Surveys, {\bf 32} (1977), 67--120.


\bibitem
{Fr} A.~Friedman, {\rm Partial Differential Equations}, Robert E.
Krieger Publ. Comp., Malabar, 1983.


\bibitem
 {Gal2m}
  V.A.~Galaktionov, {\em On
  higher-order semilinear parabolic equations with measures as initial data}, J.
  Eur. Math. Soc., {\bf 6} (2004), 193--205.
 % ({\tt http:// www.maths.bath.ac.uk/ MATHEMATICS/ preprints.html}).



%%\bibitem % [Gal2m]
%%{Gal2m1} V.A. Galaktionov, {\em On a spectrum of blow-up patterns
%%for a higher-order semilinear parabolic equations}, Proc. Royal
%%Soc. London A, {\bf 457} (2001), 1--21.

%%\bibitem
%% {GalCr}
%%  V.A. Galaktionov, {\em Critical global asymptotics in
%%  higher-order semilinear parabolic equations}, Int.
%%  J. Math. Math. Sci., {\bf 60} (2003),  3809--3825.
  %%%({\tt http://www.maths.bath.ac.uk/MATHEMATICS/
 %%%%preprints.html}).


%\bibitem
% {GH}
% V.A. Galaktinov abd P. Harwin, {\em Spectra of critical exponents
% in nonlinear heat equtions with absorption}, in preparation
% ({\tt http://www.maths.bath.ac.uk/MATHEMATICS/preprints.html}).


%\bibitem
%{GKS0}
% V.A. Galaktionov, S.P. Kurdyumov, and
%  A.A.  Samarskii, {\em On asymptotic stability of self-similar solutions
%  of the heat equation with a nonlinear sink}, Soviet Math. Dokl.,
%  {\bf 31} (1985), 271-276.



\bibitem %[GKS]
{GKS} V.A.~Galaktionov, S.P.~Kurdyumov, and
  A.A.~Samarskii, {\em On asymptotic ``eigenfunctions" of the Cauchy
  problem for a nonlinear parabolic equation}, Math. USSR Sbornik,
  {\bf 54} (1986), 421--455.



\bibitem %[GS1]
{GS1} V.A.~Galaktionov and A.E.~Shishkov, {\em Saint-Venant's
principle in blow-up for higher-order quasilinear parabolic
 equations}, Proc. Roy. Soc.
 Edinburgh, {\bf 133A} (2003), 1075--1119.


\bibitem %[GS1]
{GSVSS1} V.A.~Galaktionov and A.E.~Shishkov, {\em  Higher-order
quasilinear parabolic equations
 with singular initial data}, Comm. Contemp. Math., {\bf 8} (2006), 331--354.





%\bibitem %[GP1]
%{GP1} V.A. Galaktionov and S.I. Pohozaev, {\em Existence and
%blow-up for higher-order semilinear parabolic equations: majorizing
%order-preserving operators},
% Indiana Univ. Math. J., to appear
% ({\tt http:// www.maths.bath.ac.uk/MATHEMATICS/preprints.html}).



%%\bibitem % [GV]
%%{GV}   V.A. Galaktionov and J.L. Vazquez, {\em Continuation of
%%blow-up solutions of nonlinear heat equations in several space
%%dimensions,} { Comm. Pure Appl. Math.,} {\bf 50} (1997), 1-68.


\bibitem
{AMGV}
 V.A.~Galaktionov and J.L.~Vazquez, {\rm A Stability Technique  for Evolution Partial Differential Equations.
  A Dynamical Systems Approach},
Birkh\"auser, Boston/Berlin, 2004.

\bibitem{GW2}
V.A.~Galaktionov and J.F.~Williams,
% \bysame,
{\em On very singular similarity solutions of a higher-order
  semilinear parabolic equation}, Nonlinearity, {\bf 17} (2004), 1075--1099.
  %({\tt
  %http://www.maths.bath.ac.uk/MATHEMATICS/preprints.html}).



\bibitem
 {GohKr}
  I.C. Gokhberg and M.G. Krein, {\rm Introduction to the Theory of
  Linear Nonselfadjoint Operators}, Transl. Math. Monogr., Vol.
  {\bf 18}, AMS, Providence, RI, 1969.


%%\bibitem  %[H]
%%{He}   D. Henry, {\rm Geometric Theory of Semilinear Parabolic
%%Equations,} Lecture Notes in Math., Vol. {\bf 840},
%%Springer-Verlag, New York, 1981.


%%\bibitem  %[H]
%%{H}   D.B. Henry, {\em Some infinite-dimensional Morse-Smale
%%systems defined by parabolic partial differential equations,} { J.
%%Differ. Equat.,} {\bf 59}  (1985), 165--205.


\bibitem
 {KPApp}
 S.~Kamin and L.A.~Peletier, {\em Singular solutions of the heat
 equation with absorption}, Proc. Amer. Math. Soc., {\bf 95}
 (1985), 205--210.


\bibitem{KPe}
 S.~Kamin and L.A.~Peletier, {\em Large time behaviour of
  solutions of the porous media equation with absorption}, Israel J.
  Math., {\bf 55} (1986), 129--146.

\bibitem %[KVer]
{KVer} S.~Kamin and L.~V\'{e}ron, {\em Existence and uniqueness of
the
  very singular solution of the porous media equation with absorption},
  J. Analyse Math., {\bf 51} (1988), 245--258.


\bibitem %[Kras]
{Kras} M.A. Krasnosel'skii, {\rm Topological Methods in the Theory
of Nonlinear Integral Equations}, Pergamon Press, Oxford/Paris, 1964.



\bibitem
{KrasZ} M.A.~Krasnosel'skii and P.P.~Zabreiko, {\rm Geometrical Methods
of Nonlinear Analysis}, Springer-Verlag, Berlin/Tokyo, 1984.

 %%\bibitem
 %% {JLi}
 %% J.L. Lions, {\rm Quelques m\'{e}thodes de r\'{e}solution
%%des probl\`{e}mes aux limites non lin\'{e}aires\/}, Dunod,
%%Gauthier--Villars, Paris, 1969.



%%\bibitem
%%{Lun} A. Lunardi, {\rm Analytic Semigroups and Optimal Regularity
%%in Parabolic Problems}, Birkh\"auser, Basel/Berlin, 1995.


\bibitem
{MV02}
 M.~Marcus and L.~Veron, {\em Initial trace of positive solutions to semilinear parabolic inequalities},
{Adv. Nonlinear Studies}, {\bf 2} (2002), 395--436.


\bibitem
 {Maz}
 V.G.~Maz'ja, {\rm Sobolev Spaces}, Springer-Verlag, Berlin/Tokyo,
 1985.



%%\bibitem
%%{PelTroy} L.A. Peletier and W.C. Troy, {\rm Spatial Patterns:
%%Higher Order Models in Physics and Mechanics}, Birkh\"auser,
%%Boston/Berlin, 2001.

\bibitem
{Poh00} S.I.~Pohozaev, {\em Eigenfunctions of the equation $\Delta
u + \lambda f(u) = 0$},
 Soviet Math. Dokl., {\bf 6} (1965), 1408--1411.



\bibitem
{Rynn1}
 B.~Rynn, {\em Global bifurcation for $2m$th-order boundary
value problems and infinitely many solutions of superlinear
problems},  J. Differ. Equat., {\bf 188} (2003), 461--472.




%%\bibitem %[SGK]
%%{SGK}  A.A. Samarskii, V.A. Galaktionov, S.P. Kurdyumov, and A.P.
%%Mikhailov, {Blow-up in Quasilinear Parabolic Equations}, \rm
%%Walter de Gruyter, Berlin/New York, 1995.

\bibitem
{SW} J.~Shi and J.~Wang, {\em Morse indices and exact multiplicity
of solutions to semilinear elliptic problems}, Proc. Amer. Math.
Soc.,  {\bf 127} (1999), 3685--3695.




\bibitem %[GS1]
{SV1} A.E.~Shishkov and L.~V\'eron, 
{\em Singular solutions of some nonlinear parabolic equations with spatially inhomogeneous absorption}, Calc. Var. PDEs, {\bf 33} (2008), 343--375.
%% 35K55 (35K15)  
%%{\em Diffusion versus absorption in semilinear
%%parabolic problems}, Comptes Rendus Acad. Sci. Paris., to appear.
 %%%%%%%%, {\bf 133A} (2003), 1075--1119.


\bibitem %[Tay]
{Tay} M.E. Taylor, {\rm Partial Differential Equations III.
Nonlinear Equations,} Springer, New York/Tokyo, 1996.


\bibitem %[VainbergTr]
{VainbergTr}
 M.A.~Vainberg and V.A.~Trenogin, {\rm Theory of
Branching of Solutions of Non-Linear Equations}, Noordhoff Int.
Publ., Leiden, 1974.

%%\bibitem %[Zel]
%%{Zel} T.I. Zelenyak, {\em Stabilization of solutions of boundary
%%value problems for a second order parabolic equation with one
%%space variable}, Differ. Equat., {\bf 4} (1968), 17-22.




\end{thebibliography}
\end{document}